\documentclass[11pt]{amsart}
\numberwithin{equation}{section}
\newtheorem{theorem}{Theorem}
\newtheorem{definition}[theorem]{Definition}

\newtheorem{lemma}[theorem]{Lemma}
\newtheorem{corollary}[theorem]{Corollary}
\numberwithin{theorem}{section}

\newtheorem{remark}[theorem]{Remark}

\def\al{\aligned}
\def\eal{\endaligned}
\def\M{{\bf M}}
\def\be{\begin{equation}}
\def\ee{\end{equation}}
\def\lab{\label}
\def\a{\alpha}
\def\b{\beta}
\def\e{\epsilon}
\def\t{\tilde}

\def\M{{\bf M}}

\def\al{\aligned}

\def\g{\bar}
\def\p{\partial}
\def\d{\nabla}
\numberwithin{equation}{section}

\begin{document}

\tracingpages 1
\title[convergence]{\bf A compactness result for Fano manifolds and K\"ahler Ricci flows
}
\author{Gang Tian and Qi S. Zhang}
\address{BICMR, Peking University, Beijing, 100871, China and Department of Mathematics, Princeton University, Princeton, NJ 02139, USA}
\address{
Department of
Mathematics, University of California, Riverside, CA 92521, USA}
\date{2014/04/14}

\begin{abstract}
We obtain a compactness result for Fano manifolds and K\"ahler Ricci flows.
Comparing to the more general Riemannian versions  in
Anderson \cite{An:1} and Hamilton \cite{Ha:1}, in this Fano case,
the curvature assumption is much weaker and is preserved by the
K\"ahler Ricci flows. One assumption
 is the boundedness of
the Ricci potential and the other is the smallness of Perelman's
entropy. As one application, we obtain a new local regularity
criteria and structure result for K\"ahler Ricci flows. The proof
is based on a H\"older estimate for the gradient of harmonic
functions, which may be of independent interest.

\end{abstract}
\maketitle
\tableofcontents
\section{Introduction}

Compactness theorems have been useful tools in the study of
geometric objects such manifolds and their evolutions under
certain equations. Well known examples include the Cheeger-Gromov
compactness theorem in the case of bounded curvature, diameter and
volume lower
 bound (\cite{Ch:1}, \cite{G:1}), as well as M. Anderson's extension to the
case of bounded Ricci curvature, diameter and injectivity lower
bound (\cite{An:1}). In both cases, the compactness is in the
$C^{1, \a}$ topology for any $\a$ in $(0, 1)$. M. Anderson also
mentioned that if one replaces the $L^\infty$ bound of the Ricci
curvature by its $L^p$ bound with $p>n/2$, then the compactness
holds in $C^\a$ topology for some $\a \in (0, 1)$. See also
\cite{PW:1}

In this paper we prove a similar compactness result for Fano
manifolds. It is shown that in this special case, the assumption
on the Ricci curvature can be much weaker and preserved by the
normalized Ricci flow. Then we will apply this result to study
K\"ahler Ricci flows.

In order to present the results precisely, We will use the
following notations and definitions. We use $\M$ to denote a $n$
real dimensional compact Fano manifold, i.e. K\"ahler manifold
with positive first Chern class. Denote by $g$ be a stationary
metric and by $g(t)$ the metric at time $t$; $d(x, y, t)$ is the
geodesic distance under $g(t)$; $B(x, r, t) = \{ y \in {\M} \ | \
d(x, y, t) < r \}$ is the geodesic ball of radius $r$, under
metric $g(t)$, centered at $x$, and $|B(x, r, t)|_{g(t)}$ is the
volume of $B(x, r, t)$ under $g(t)$;  $d g(t)$ is the volume
element. We also reserve $R=R(x, t)$ as the scalar curvature under
$g(t)$. When the time variable $t$ is not explicitly used, we may
also suppress it in the notations mentioned above. In this paper,
when mentioning a single manifold,
 it is always an $n$ real dimensional compact Fano manifold $\M$ which satisfies
  the following

\noindent {\bf Basic assumptions}:

{\it Assumption  1. $L^2$ Sobolev inequality: there is a positive constant $s_0$ such that
\[
 \left( \int_{\M}  v^{2n/(n-2)} d g \right)^{(n-2)/n} \le s_0 \left( \int_{\M} | \nabla v |^2
 dg +  \int_{\M} v^2
 dg \right)
\]for all $ v \in C^\infty(\M)$.
}

{\it Assumption  2. There exist positive constants $\kappa$ and $D$, such that
\[
\kappa r^n \le |B(x, r)| \le \kappa^{-1} r^n, \qquad \forall x \in \M, \quad 0 < r< diam(\M) \le D.
\]}

{\it Assumption 3. Let $u$ be the Ricci potential i.e. $R_{i \g j}
= - u_{i \g j} +g_{i \g j}$. Then $u$ satisfies $\Vert u
\Vert_{C^1(\M)} \le E.$ Moreover the scalar curvature satisfies $\Vert R \Vert_\infty \le E$.}

One motivation for the basic assumptions is that they are
satisfied, with uniform constant, by metrics along a normalized K\"ahler Ricci flow. See
Property KRF below. Also the volume lower bound is redundant since the
Sobolev inequality implies it.

Parts of our results are also related to Perelman's $W$ entropy by

\begin{definition}
\lab{defWshang} Let $D$ be a domain in $\M$. The Perelman $W$ entropy (functional) with
parameter $\tau>0$ is the quantity
\be
W(g, f, \tau) = \int_D [\tau (R + |\d f|^2) + f - n] e^{-f}  (4 \pi \tau )^{-n/2} dg, \quad
f \in C^\infty_c(D).
\ee The $\mu$ invariant with parameter $\tau$ is the infimum of the $W$ entropy,
given by
\be
\mu(g, \tau, D) = \inf \{ W(g, f, \tau) \, | \, f \in C^\infty_c(D), \,
(4 \pi \tau )^{-n/2} \int_D e^{-f} dg =1 \}.
\ee
\end{definition}

It is clear that $\mu(g, \tau, D) \ge \mu(g, \tau, \M )$. According to Perelman
\cite{P:1}, the latter is monotone non-decreasing under Ricci flow if the
parameter $\tau$ satisfies $\frac{d \tau}{d t} = -1$ . If no confusion arises, we
will drop the argument $\M$ and simply write $\mu(g, \tau, \M )$ as $\mu(g, \tau)$.

\medskip

The first main results of the paper is the following theorem.

\begin{theorem} (compactness of Fano manifolds)
\lab{thmconv} (a). Let $SP=SP(s_0, \kappa, D, E, i_0)$ be the space of compact Fano manifolds $\M$ satisfying Assumptions 1, 2, 3
 and that $inj_{\M} \ge i_0$. Then there exists a positive number $\a \in (0, 1)$ such that
 $SP$ is compact
in the $C^{1, \a}$ topology.

(b).  Let $SP=SP(s_0, \kappa, D, E, \eta_0, \tau_0, r_0)$ be the space of compact Fano manifolds
 $({\M}, g)$
satisfying Assumptions 1, 2, 3 and that $\sup_{\tau \in (0, \tau_0], \, p \in \M}\mu(g, \tau, B(p, r_0)) \ge - \eta_0$.
Here $\eta_0, \tau_0, r_0$ are any given positive numbers.
If $\eta_0$ is sufficiently small, then there exists a positive number $\a \in (0, 1)$ such that
 $SP$ is compact
in the $C^{1, \a}$ topology.
\end{theorem}

Next, we apply the previous theorem to the
 study  the (normalized) K\"ahler Ricci flows
\be
\lab{krf}
\partial_t g_{i\g j} = -  R_{i\g j} + g_{i\g{j}} = \p_i \p_{\g j} u, \quad t>0,
\ee on a compact, K\"ahler manifold $\M$ of complex dimension
$m=n/2$, with positive first Chern class. We always assume that
the initial metric is in the canonical K\"ahler class $2\pi
c_1(\M)$.

Given initial K\"ahler metric $g_{i\g j}(0)$, H. D. Cao
\cite{Ca:1} proved that (\ref{krf}) has a solution for all time
$t$. Recently, many  results concerning long time and uniform
behavior of (\ref{krf}) have appeared. For example, when the
curvature operator or the bisectional curvature is nonnegative, it
is known that solutions to (\ref{krf}) stays smooth when time goes
to infinity (see \cite{CCZ:1}, \cite{CT:1} and \cite{CT:2} for
examples). In the general case, Perelman (cf \cite{ST:1}) proved
that the scalar curvature $R$ is uniformly bounded, and the Ricci
potential $u(\cdot, t)$ is uniformly bounded in $C^1$ norm, with
respect to $g(t)$. When the complex dimension $m=2$, let $({\M},
g(t))$ be a solution to (\ref{krf}),  it is proved in
(\cite{CW:1}) that the flow sequentially converges to an orbifold.
When $m=3$, the authors of \cite{TZz:1} proved sequential
convergence except on a small singular set. For the general case,
in the paper \cite{TZq:1}, we proved that a Gromov-Hausdorff limit
of the flow is a metric space with volume doubling
 and $L^2$ Poincar\'e inequality and $L^1$ isoperimetric inequality.
 So it is a PI space in the sense of Cheeger.
There are many more interesting papers  on the convergence issue
of K\"ahler Ricci flow. See for example, \cite{CLW:1},
\cite{MS:1}, \cite{PS:1}, \cite{PSSW1:1}, \cite{PSSW2:1},
\cite{Se:1}, \cite{SoT:1}, \cite{Sz:1}, \cite{SW:1}, \cite{To:1},
\cite{TZhu1:1}, \cite{TZhu2:1}, \cite{TZZZ:1}, \cite{Zz1:1},
\cite{Zz2:1} and \cite{Zhu:1} and references therein.

The following is a local regularity criteria for K\"ahler Ricci flows,
which can be regarded
as a strengthened pseudolocality theorem for K\"ahler Ricci flow.
Note our condition on the Ricci curvature is preserved under the flow.
Hence essentially
 the result does not require any direct assumption on Ricci curvature or
curvature tensor except on the initial metric. Some applications will be given in the form of
 corollaries, including the compactness result.

\begin{theorem} (local regularity criteria for K\"ahler Ricci
flows) \lab{thSLocKRF} Let $({\M, g(t)})$ be a normalized K\"ahler
Ricci flow on a compact Fano manifold whose initial metric
satisfies the basic assumptions. There exists positive
 numbers $r_0$ and $\eta$ such that the following statement is true.

Suppose either one of the conditions holds for $r \in (0, r_0]$.

1. the geodesic ball $B(x, r, 0)$ is almost Euclidean, i.e.  for a sufficiently
small positive number $\eta$, any $y \in B(x, r, 0)$ with $B(y, \rho, 0) \subset B(x, r, 0)$, it holds
\be
|B(y, \rho, 0)|_{g(0)} \ge (1-\eta) \omega_n \rho^n
\ee where $\omega_n$ is the volume of the unit ball in ${\bf R}^n$;

2. for a sufficiently small positive number $\eta>0$ and a given number $\tau_0>0$,
 the infimum of the  $W$ entropy on $B(x,  r, 0)$
satisfies \be \sup_{\tau \in (0, \tau_0]} \mu(g(0), \tau, B(x, r,
0 )) \ge -\eta. \ee

3. the injectivity radius of the ball $B(x, r, 0)$ is bounded from below by a number $i_0>0$.

Then, exists a number $\e>0$ such that
\be
|Rm(y, t)| \le t^{-1} + (\e r)^{-2}
\ee for all $t \in (0, (\e r)^2]$ and $y$ such that $d(x, y, t) \le \e r.$
\end{theorem}

\medskip

\begin{remark}
(a). The result under condition 1 is an improvement to Theorem 4.2
in \cite{TZz:2} where the same conclusion is reached under the
extra assumption that $|\d \d u|$ is $L^\infty$ uniformly in time.

(b). The result under condition 2  is similar in spirit to the new
$\e$ regularity result for Ricci flow in \cite{HN:1} where the
authors proved boundedness of the curvature tensor in a space-time
cube under the condition that a heat kernel weighted entropy is
close to zero. Since the heat kernel in that result is coupled to
the Ricci flow, that condition is a space time one. In contrast
our condition is applied on the initial value.

(c). The result under condition 2 is also related to the
Pseudo-locality Theorem 3.1 in \cite{TW:1} where a Ricci lower
bound is assumed for general normalized Ricci flows.

(d). The result also holds for un-normalized Ricci flows by
scaling and choosing $\e$ sufficiently small.

(e). One may prove the result under condition 3 by different
methods.
\end{remark}

One application of the theorem is the following compactness result for
K\"ahler Ricci flows.

\begin{corollary} (compactness of  K\"ahler Ricci
flows) \lab{cocompactKRF} Let $({\M}_k, g_k(t), p_k)$ be a
sequence of normalized compact K\"ahler Ricci flows in the time
interval $[0, T]$. Then the following conclusion holds.

(a). Suppose the initial metrics $({\M}_k, g_k(0))$ satisfies the basic assumptions and
 that the injectivity radius of $({\M}_k, g_k(0))$ is bounded from below by a positive constant
$i_0$. Then, there exists a positive constant $T_1 \le T$ such that $({\M}_k, g_k(t), p_k)$,
$t \in (0, T_1]$,
is compact in $C^\infty_{loc}$ topology.

(b).  Suppose the initial metrics $({\M}_k, g_k(0))$ satisfies the basic assumptions.
Let $\{ t_j \}$ be a partition of the time interval such that $\sup (t_j-t_{j-1}) \le T_1$.
 Suppose the injectivity radii of $({\M}_k, g_k(t_j))$ are bounded from below by a positive constant
$i_0$. Then $({\M}_k, g_k(t), p_k)$,
$t \in (0, T]$,
is compact in $C^\infty_{loc}$ topology.
\end{corollary}

\begin{remark}
Recall that Hamilton's  compactness result for the general Ricci
flow in forward time requires uniform bound for the curvature
tensor through out space time and the positive lower bound of the
injectivity radii at initial time. By the classical volume
comparison theorem of Gunther, the initial manifolds are $\kappa$
non-collapsed. Now, by Perelman \cite{P:1}, we know that the Ricci
flows are also $\kappa$ non-collapsed in the time interval $[0,
T]$. Hence the injectivity radius are bounded from below in all
space time by the bound on the curvature tensor. Therefore,
Hamilton's compactness theorem for general Ricci flows implicitly
requires both curvature bound and injectivity lower bound in all
space time.

In contrast, part (a) of the corollary for K\"ahler Ricci flow only requires lower bound of injectivity radius
and the basic assumptions on the initial metric. Part (b) requires lower bound of
injectivity radius
for some discrete times
and the basic assumptions on the initial metric. So, in the K\"ahler Ricci flow case,
we have essentially removed the direct assumption
on curvature tensor from compactness result.
\end{remark}

With possible application in mind, we single out another
consequence of Theorem \ref{thSLocKRF} as

\begin{corollary}  ($C^{1, \a}$  structure of almost Euclidean region).
Let $({\M}, g(t))$ be a normalized K\"ahler Ricci flow whose
initial metric $g(0)$ satisfies the basic assumptions. There exist
positive numbers $\eta \in (0, 1)$, $\delta$ and $r_0$ such that
the following
 statement holds.

For any $t>0$, $x \in \M$ and $r \in (0, r_0]$, suppose the ball $B(x,
r, g(t))$ is almost Euclidean in volume, i.e. for all $B(y, \rho,
g(t)) \subset B(x, r, g(t))$, \be |B(y, \rho, g(t))|_{g(t)} \ge
(1-\eta) \omega_n \rho^n,
\ee where $w_n$
is the volume of $n$ dimensional Euclidean unit ball.
 Then,
\[ r^{\a, \theta}_{g(t)}(x) \ge
\delta r
\]Here $r^{\a, \theta}_{g(t)}(x)$ is the $C^{1, \a}$ harmonic radius
at $x$ defined in Definition \ref{defharmrad} in Section 4. The $C^{1, \a}$ norm
of the metric is  bounded within the radius.
\end{corollary}

The proof follows from part 2 of Lemma \ref{leharmrad} which is
the core of Theorem \ref{thSLocKRF} and the fact that the basic
assumptions are preserved under the K\"ahler Ricci flow. See  the
Property KRF below.

\begin{remark}
 We mention that a similar $C^\a$ structure result is given by
Theorem 2.35 in \cite{TZz:1}, under the assumption that the $L^p$, $p>n/2$ norm
of the Ricci curvature is uniformly bounded. That result plays an important role in
their proof of convergence result in complex dimension 3.
\end{remark}

The proof of Theorem \ref{thSLocKRF}  is based on the following properties for  K\"ahler Ricci flow on a compact manifold with positive first Chern class.

\vskip 0.5cm

\noindent {\it Property KRF.  Let $({\M}, g(t))$ be  a K\"ahler Ricci flow (\ref{krf}) on a compact manifold with positive first Chern class.
There exist uniform positive constants $C$ and $\kappa$ depending only on
$g(0)$ so that \\

1.  $| R(g(t))| \le C,$

2. $diam (\M, g(t)) \le C,$

3. $\Vert u(\cdot, t) \Vert_{C^1} \le C.$

4. $|B(x, r, t)|_{g(t)} \ge \kappa r^n$, for all $t>0$ and $r \in (0, diam (\M, g(t)))$. }

5. $
|B(x, r, t)|_{g(t)} \le \kappa^{-1}  r^n
$ for all $r>0$, $t>0$.

6. There exists a uniform constant $S_2$ so that the following $L^2$ Sobolev inequality holds: \\
 \[
 \left( \int_{\M}  v^{2n/(n-2)} d g(t) \right)^{(n-2)/n} \le S_2 \left( \int_{\M} | \nabla v |^2
 dg(t) +  \int_{\M} v^2
 dg(t) \right)
\]for all $ v \in C^\infty(\M, g(t))$.

7. (a). Let $\Gamma$ be the Green's function on $({\M}, g(t))$.
Then there exists a uniform constant $C$ such that \be
\lab{GdGjie} |\Gamma(x, y) | \le \frac{C}{d(x, y)^{n-2}}, \quad
|\nabla \Gamma(x, y) | \le \frac{C}{d(x, y)^{n-1}}. \ee

(b). Let $p=p(x, y, s)$ be the (stationary) heat kernel for
$({\M}, g(t))$.  There exist positive constants $a_1$ and $a_2$,
depending only on $g(0)$ such that \be \frac{a_1}{s^{n/2}} e^{-a_2
d(x, y)^2/s} \le p(x, y, s) \le \frac{1}{a_1 s^{n/2}} e^{- d(x,
y)^2/(a_2 s)}, \qquad s \in (0, 1]; \ee

 \be |\d p(x, y, s)| \le \frac{1}{a_1 s^{(n+1)/2}}
e^{- d(x, y)^2/(a_2 s)}, \qquad s \in (0, 1]. \ee

8. uniform $L^2$ Poincar\'e inequality: for any $v \in
C^\infty(B(x, r))$ where $B(x, r)$ is a proper ball in $({\M},
g(t))$, there is a uniform constant $C$ such that \be \lab{l2PI}
\int_{B(x, r)} |v-  v_B |^2 dg \le C r^2 \int_{B(x, r)} |\d v|^2 dg.
\ee Here all quantities are with respect to  $g=g(t)$ and $v_B$
is the average of $v$ in $B(x, r)$.

Property A  1-4 is due to Perelman (c.f. \cite{ST:1}); Property 5 can be found in \cite{Z11:1}
and also \cite{CW:2}; Property 6 was first proven in
\cite{Z07:1} (see also \cite{Ye:1}, \cite{Z10:1} ).
 Properties 7 and 8 are in \cite{TZq:1}. Bounds on the Green's function and
heat kernel are well known if the metrics have uniform $L^\infty$
 lower bound for the Ricci curvature (\cite{LY:1}),  a condition that is unavailable here.

The rest of the paper is organized as follows. In Section 2, we will prove some
integral bounds for the Hessian of the Ricci potential, which implies that
the Ricci curvature is actually small in certain Morrey or Kato type norm.
In Section 3, we show that the gradient of the Ricci potential is H\"older
continuous within a harmonic coordinate. In Section 4, we will prove a
lower bound for the harmonic radius under three separate conditions. One
involves the injectivity radius, another involves volume of balls and third one
relates with Perelman's W entropy. The theorems and corollaries will be proven in Section 5.
Sometimes we need to switch between real and complex coordinates in computations, which
may result in extra harmless constants in the Laplacian, Ricci curvature, etc.
\medskip

\section{Some integral estimate on the Hessian of the Ricci potential}

In order to prove the theorems, in this section we state and prove
an a priori integral estimate for the Hessian of the Ricci
potential, which essentially means that it is a sub-critical
quantity  comparing with the Laplacian. In fact this  simple
estimate already shows that in integral sense $| Ric |$ scales
like $1/r$ instead of the usual $1/r^2$. Here $r$ is the distance
function. It is the reason that one gains one order of regularity
in the Theorems.

\begin{lemma}
\lab{lehess}
Let $u$ be the Ricci potential. Let $B(x, r)$ be a proper geodesic ball with $r \le 1$.
Then there exists a uniform constant $A_0$ depending only on the
parameters in the basic assumptions of the manifold and $\Vert \nabla u \Vert_\infty$ and $\Vert R \Vert_\infty$
such that
\[
\int_{B(x, r)} | Hess \, u |^2 dg \le A_0 r^{n-2}.
\]
\end{lemma}
\proof Starting  from the equation, in real variable form
 \be
\Delta u = n - R, \ee and applying the Bochner's formula, we know
that \be \lab{bochn} \Delta | \nabla u |^2 = 2 | Hess u |^2 + 2
R_{i \g j} u_j u_{\g i} - 2 \nabla R \nabla u \ee Here $\Delta $
is the real Laplacian. Let  $\phi$ be a standard Lipschtiz cut-off
function such that $\phi=1$ on $B(x, r)$, $\phi=0$ on $B(x, 2r)^c$
and $|\nabla \phi| \le C/r$. After integration, we have \be
\lab{hess1} \al
&2 \int | Hess \, u |^2 \phi^2 dg\\
 &=\int \Delta | \d u|^2 \, \phi^2 dg - 2 \int
R_{i \g j} u_j u_{\g i} \phi^2 dg + \int 2 \d R \d u \phi^2 dg\\
&=-4 \int < \d | \d u |, \d \phi > |\d u | \phi dg + 2 \int (\partial_i \partial_{\g j}
u - g_{i \g j} ) u_j u_{\g i} \phi^2 dg \\
&\qquad - 2 \int R (n-R) \phi^2 dg -
4 \int R \d u \d \phi \phi dg.
\eal
\ee Note that $| \d | \d u | | \le | Hess \, u |$ and,  by Perelman (cf \cite{ST:1}), that
$|\d u|$ and $|R|$ are bounded. We can apply Cauchy-Schwarz inequality to deduce
\be
\lab{hess2}
\al
& \int | Hess \, u |^2 \phi^2 dg \\
&\le \frac{C}{r^2} \int_{B(x, 2r)} (|\d u|^2 +
|\d u|^4) dg + 2 \int_{B(x, 2r)}  |R (n-R)| dg +
 \frac{C}{r}  \int_{B(x, 2r)} |R| |\d u| dg \le C r^{n-2}.
\eal
\ee To get the last inequality, we have used the $\kappa$ non-inflating property.
\qed

The next lemma provides a lower bound on the Green's function $\Gamma=
\Gamma(x, y)$. Since $\M$ is compact, we know that $\Gamma$ changes sign.
So the lower bound holds only when $x$ and $y$ are close.

\begin{lemma}
\lab{leGlow}
Let $\Gamma$ be the Green's function of the Laplace operator on scalar functions.
There exist positive numbers $r_0$ and $C$, depending only on the parameters of the
basic assumptions on $\M$, such that
\be
\Gamma(x, y) \ge \frac{C}{d(x, y)^{n-2}}
\ee provided that $d(x, y) \le r_0$.
\proof
\end{lemma}

Since $\M$ is a compact manifold, it is well known that \be
\lab{gamma-G} \Gamma(x, y) = \int^\infty_0 \left( G(x, t, y) -
\frac{1}{|\M|} \right) dt \ee where $G$ is the heat kernel on
$\M$. We mention that the time $t$ here is not the time in the
Ricci flow. From Section 2 in \cite{TZq:1}, the following
inequalities hold for $G$. There exist positive constants $c_1,
..., c_4$ and $\b$ such that \be \left| \int^\infty_\b [ G(x, t,
y) - \frac{1}{|\M|} ] dt \right| \le c_1 \int^\infty_{\b} e^{- c_2
t} dt \le c_1/c_2; \ee Also, if $t \in (0, \b)$, then \be G(x, t,
y) \ge \frac{c_3}{t^{n/2}} e^{-c_4 d^2(x, y)/t}. \ee Substituting
these into (\ref{gamma-G}), we find that \be \al \Gamma(x, y) &\ge
\int^\b_0 \frac{c_3}{t^{n/2}} e^{-c_4 d^2(x, y)/t} dt
- \b/|\M| - c_1/c_2\\
&=\frac{c_3}{d(x, y)^{n-2}} \int^{\b/d(x, y)^2}_0 \frac{1}{s^{n/2}}
e^{-c_4/s} ds - \b/|\M| - c_1/c_2.
\eal
\ee Hence, there exists $r_0>0$ such that
\[
\Gamma(x, y) \ge \frac{C}{d(x, y)^{n-2}}
\] provided that $d(x, y) \le r_0$. \qed

\begin{lemma}
\lab{lehess2}
Let $u$ be the Ricci potential.
Then there exists a uniform constant $A_1$ depending only on the
parameters of the basic assumptions of the manifold such that
\be
\lab{Khess}
K(|Hess \, u|^2) \equiv \sup_{x \in \M} \int_{\M}
\frac{ | Hess \, u(y) |^2}{d(x, y)^{n-2}} dg(y) \le A_1.
\ee Here and later we have used the notation for the Kato norm of an integrable
function $F$,
\be
K(F)=\sup_{x \in \M} \int_{\M}
\frac{ |F(y)|}{d(x, y)^{n-2}} dg(y).
\ee
\proof
\end{lemma}
Recall from (\ref{bochn}) that
\[
\Delta | \nabla u |^2 = 2 | Hess u |^2 + 2 R_{i \g j} u_j u_{\g i} - 2 \nabla R \nabla u.
\] Let $a$ be the average of $|\d u|^2$ over $\M$.  Then, from the definition of
Green's function, we have
\be
\al
| \nabla u |^2(x) - a &= - 2 \int \Gamma(x, y) | Hess u (y)|^2 dg(y) - 2 \int
\Gamma(x, y) R_{i \g j} u_j u_{\g i} (y) dg(y) \\
&\hskip2cm+
2 \int \Gamma(x, y) \nabla R \nabla u (y) dg(y).
\eal
\ee Using the relation $\p_i\p_{\g j} u = g_{i\g j}-R_{i\g j}$ and applying integration by parts, we see that
\be
\lab{Ghess}
\al
2 &\int \Gamma(x, y) | Hess u (y)|^2 dg(y) \\
&=
a- | \nabla u |^2(x) + 2 \int \Gamma(x, y) (\p_i\p_{\g j} u - g_{i\g j}) u_j u_{\g i} (y) dg(y)
+
2 \int \Gamma(x, y) \nabla R \nabla u (y) dg(y)\\
&=(a- | \nabla u |^2(x)) - 2 \int u_i \Delta u  u_{\g i} \Gamma(x, y) dy
- 2 \int u_i u_j u_{\g i \g j} \Gamma(x, y) dg(y) \\
&\qquad
-2 \int u_i u_j u_{\g i} \d_{\g j} \Gamma(x, y) dg(y)
- 2 \int |\d u|^2 \Gamma(x, y) dg(y) \\
&\qquad - 2 \int R \Delta u \, \Gamma(x, y) dg(y)
-2 \int R \d u \d \Gamma(x, y) dg(y)\\
&\equiv I_1 + ...+ I_7.
\eal
\ee Our next task is to bound each term $I_k$, $k=1, ..., 7$.

According to Perelman $|\d u|$ and $R$ are bounded. Hence $|I_1|
\le C$ and \be \lab{I2} |I_2 | \le 2 \int |\d u | \, | n-R| \,
|\Gamma(x, y) dg(y) \le C \int \frac{1}{d(x, y)^{n-2}} dg(y). \ee
Here we just used the bound $|\Gamma(x, y) |  \le \frac{C}{d(x,
y)^{n-2}}$, which was proven in \cite{TZq:1}. From (\ref{I2}) and
the volume non-inflating property, it is easy to see that $|I_2|
\le C$. Therefore, \be \lab{I12} |I_1| + |I_2 | \le C. \ee Using
the bound on $\Gamma(x, y)$ again, we have
 \be \lab{I3} |I_3| \le
C \int | \d u |^2 \, | Hess \, u | \frac{1}{d(x, y)^{n-2}} dg(y) \le C.
\ee

Applying the gradient bound  $|\d \Gamma(x, y) |  \le \frac{C}{d(x, y)^{n-1}}$
 in \cite{TZq:1}, we find
\be
\lab{I4}
|I_4| \le C \int \frac{|\d u |^3}{d(x, y)^{n-1}} dg(y) \le C
\ee where we have used the volume non-inflating property again.

Similarly,
\be
\lab{I567}
|I_5| +|I_6| + |I_7| \le C + C \int \frac{ | R (n-R)(y) |}{ d(x, y)^{n-2}}
+ C \int \frac{|R \, \d u |}{d(x, y)^{n-1}} dg(y) \le C.
\ee Substituting (\ref{I12}), (\ref{I3}), (\ref{I4}) and (\ref{I567}) into (\ref{Ghess}), we deduce
\be
\lab{Ghess2}
2 \int \Gamma(x, y) | Hess u (y)|^2 dg(y) \le  C + C \int | \d u |^2 \, | Hess \, u | \frac{1}{d(x, y)^{n-2}} dg(y).
\ee

On the other hand, by Lemma \ref{leGlow}, there are uniform constants $r_0, C>0$
such that $\Gamma(x, y) \ge C/d(x, y)^{n-2}$ when $d(x, y) \le r_0$.
Hence
\be
\al
\int_{d(x, y) \le r_0}& \frac{ | Hess \, u (y) |^2}{d(x, y)^{n-2}}
\le \frac{1}{C}\int_{d(x, y) \le r_0} | Hess \, u (y) |^2 \Gamma(x, y) dg(y) \\
&=
\frac{1}{C}\int_{\M} | Hess \, u (y) |^2 \Gamma(x, y) dg(y) -
\frac{1}{C}\int_{d(x, y) \ge r_0} | Hess \, u (y) |^2 \Gamma(x, y) dg(y)\\
&\le \frac{1}{C}\int_{\M} | Hess \, u (y) |^2 \Gamma(x, y) dg(y) +
C \int_{d(x, y) \ge r_0} | Hess \, u (y) |^2 \frac{1}{d(x,
y)^{n-2}} dg(y). \eal \ee Here we just used the upper bound
$|\Gamma(x, y) | \le C/d(x, y)^{n-2}$ in \cite{TZq:1}. This bound
and (\ref{Ghess2}) then imply that \be \lab{Gd2-n} \int_{d(x, y)
\le r_0} \frac{ | Hess \, u (y) |^2}{d(x, y)^{n-2}} \le C
\int_{d(x, y) \ge r_0} \frac{ | Hess \, u (y) |^2}{d(x, y)^{n-2}}
+ C + C \int_\M \frac{| \d u |^2 \, | Hess \, u(y) |}{d(x,
y)^{n-2}} dg(y). \ee This implies \be \lab{Gd2-nM} \al
 &\int_{\M} \frac{ | Hess \, u (y) |^2}{d(x, y)^{n-2}} \\
&\le (C+1) \int_{d(x, y) \ge r_0} \frac{ | Hess \, u (y) |^2}{d(x, y)^{n-2}} +
C + C \int_\M \frac{ | \d u |^2 \, | Hess \, u(y) | }{d(x, y)^{n-2}} dg(y)\\
&\le (C+1) r^{-(n-2)}_0 \int_{\M}  | Hess \, u (y) |^2 +
C + C \int_\M \frac{ | \d u |^2 \, | Hess \, u(y) | }{d(x, y)^{n-2}} dg(y).
\eal
\ee In Lemma \ref{lehess} we take $r=r_0$. By the uniform diameter bound
of $\M$, we know that $\int | Hess \, u |^2 dg <C$, which together with (\ref{Gd2-nM}), show that
\be
\int_{\M} \frac{ | Hess \, u (y) |^2}{d(x, y)^{n-2}} \le  C +
C \int_\M \frac{ | \d u |^2 \, | Hess \, u(y) | }{d(x, y)^{n-2}} dg(y).
\ee
After using Cauchy-Schwarz inequality, we deduce
\be
\int_{\M} \frac{ | Hess \, u (y) |^2}{d(x, y)^{n-2}} \le  C +
C \int_\M \frac{ | \d u |^4 }{d(x, y)^{n-2}} dg(y) \le C.
\ee This proves the lemma.
\qed

The next lemma is an embedding result which implies that $|Hess u|^2$, regarded
as a potential function or inhomogeneous term in an equation,
 is dominated
by the Laplacian.

\begin{lemma}
\lab{leimbed}
 Let $V$ be a smooth function on $\M$, $p$ be a point on $\M$ and $r$ be a
positive number such that $r \le diam(\M)/2$. Then for any smooth function $f$
on $\M$, the following embedding result holds
\be
\al
&\int_{B(p, r)} |V(x)| f^2(x) dg(x) \\
&\le C \sup_{z \in B(p, 2r)} \int_{B(z, 2r)}
\frac{|V(x)|}{d(z, x)^{n-2}} dg(x) \left( \Vert \d f \Vert^2_{L^2(B(p, 2r)}
+ r^{-2} \Vert  f \Vert^2_{L^2(B(p, 2r)} \right).
\eal
\ee Here $C$ is a positive constant depending only on the parameters in the basic
assumption for $\M$.

In particular, for the Ricci potential $u$, it holds
\be
\int_{B(p, r)} |Hess \, u|^2 f^2(x) dg(x)
\le C K(|Hess \, u|^2)  \left( \Vert \d f \Vert^2_{L^2(B(p, 2r)}
+ r^{-2} \Vert  f \Vert^2_{L^2(B(p, 2r)} \right).
\ee Here $ K(|Hess \, u|^2)$ is defined in (\ref{Khess}).
\proof
\end{lemma}

With the gradient bound for the Green's function and volume lower
and upper bound in hands, the lemma and its proof are essentially
known in the literature. See \cite{Si:1} and \cite{CGL:1} for
example. The only difference is that we are dealing with compact
manifolds which create one extra term from the Green's formula.
Therefore we will just sketch the proof.

Without loss of generality we assume $V \ge 0$.
Let $\phi$ be a Lipschitz cut-off function in $B(p, 2r)$ such that $\phi =1$ on
$B(p, r)$ and that $|\d \phi | \le C/r$. Then from the Green's formula, we have
\be
f\phi(x) - ave (f \phi) = - \int \Gamma(x, y) \Delta (f \phi)(y) dg(y)
\ee where $ave (f \phi)$ is the average of $f \phi$ over $\M$, and $\Gamma$ is the
Green's function. After integration by parts, this becomes
\be
f\phi(x) = ave (f \phi) + \int \d_y \Gamma(x, y) \d (f \phi)(y) dg(y).
\ee By the bound for $\d \Gamma$ in (\ref{GdGjie}), this infers
\be
|f \phi(x)| \le |ave (f \phi)| + \int \frac{C |\d f| \phi(y)}{d(x, y)^{n-1}}  dg(y)
+ \int \frac{C | f \d \phi(y|)}{d(x, y)^{n-1}}  dg(y).
\ee

For any smooth test function $\eta$ supported in $B(p, r)$, we have, from the previous inequality,
\be
\lab{Vfi123}
\al
&\int \sqrt{V} f \eta(x) dg(x)  \\
&\le C \int \sqrt{V} \eta (x) \int \frac{ |\d f| \phi(y)}{d(x, y)^{n-1}}  dg(y) dg(x) + C \int \sqrt{V} \eta (x) \int \frac{ | f \d \phi(y)|}{d(x, y)^{n-1}}  dg(y) dg(x) \\
&\qquad  +  \int \sqrt{V} \eta (x) ave ( f \phi) dg(x)\\
&\equiv I_1 + I_2 +I_3. \eal
\ee
Following the argument in
\cite{Si:1} and \cite{CGL:1}, one can show that
\be \lab{Vfi1}
I^2_1 \le C \Vert \d f \Vert^2_{L^2(B(p, 2r))} \Vert \eta
\Vert^2_{L^2(B(p, r)} \sup_{z \in B(p, 2r)} \int \frac{V(x)}{d(x,
z)^{n-2}} dg(x). \ee \be \lab{Vfi2} I^2_2 \le \frac{C}{r^2} \Vert
\d f \Vert^2_{L^2(B(p, 2r))} \Vert \eta \Vert^2_{L^2(B(p, r)}
\sup_{z \in B(p, 2r)} \int \frac{V(x)}{d(x, z)^{n-2}} dg(x). \ee
Next \be \al
|I_3| &\le \int \sqrt{V} |\eta(x)| \frac{1}{|\M|} \int |f| \phi(y) dg(y) dg(x)\\
&\le \int \sqrt{V} |\eta(x)| \frac{diam(\M)^{n-1}}{|\M|} \int \frac{|f| \phi(y) dg(y)}{d(x, y)^{n-1}} dg(x)\\
&\le C \int \sqrt{V} |\eta(x)|  \int \frac{|f| \phi(y) dg(y)}{d(x, y)^{n-1}} dg(x),
\eal
\ee where we have used Perelman's bound on the diameter and volume non-collapsing property.  Just like the case for $I_1$, we can now deduce
\be
\lab{Vfi3}
|I_3|^2 \le C \Vert  f \Vert^2_{L^2(B(p, 2r))} \Vert \eta \Vert^2_{L^2(B(p, r)}
\sup_{z \in B(p, 2r)}
\int \frac{V(x)}{d(x, z)^{n-2}} dg(x).
\ee Substituting (\ref{Vfi1}), (\ref{Vfi2}) and (\ref{Vfi3}) into (\ref{Vfi123}),
we find that
\be
\al
&\int \sqrt{V} f \eta(x) dg(x) \\
&\le
C \left[\Vert \d f \Vert_{L^2(B(p, 2r))}  + \frac{1+r}{r} \Vert  f \Vert_{L^2(B(p, 2r))} \right] \Vert \eta \Vert_{L^2(B(p, r)}
\left(\sup_{z \in B(p, 2r)}
\int \frac{V(x)}{d(x, z)^{n-2}} dg(x) \right)^{1/2}.
\eal
\ee Since $\eta$ is arbitrary, the result follows by applying the Riesz theorem. \qed

\section{$C^{1, \a}$ bounds for harmonic functions and Ricci potential in
harmonic charts}

The goal of this section is to prove $C^{1, \a}$ bounds for a
harmonic function $h$ and Ricci potentials in a harmonic
coordinate chart. We will use De Giorgi's method. In doing so, we
first need to tackle two technical issues. One is that the De
Giorgi method does not work for systems of equations in general.
This can be handled since the system for $d h$ is a weakly coupled
one involving the Ricci curvature and by the integral bound for
the Ricci curvature in the previous section. The other issue is
that we need a weighted $L^p$ Poincar\'e type inequality for some
$p$ strictly less than $2$.

\begin{lemma}
\lab{lewfPI} Given a proper geodesic ball $B=B(p, r)$ and two
functions $v \in C^\infty(B)$ and $ 0 \le f \in L^\infty(B)$, there
exist a positive number $p_0<2$ and constant $C>0$ such that
\be
\int_B |v - v_f|^{p_0} dg
\le C r^{p_0} \left[ \sup f^{p_0} \frac{|B| |B \cap supp \,
f|^{p_0-1}}{\Vert f \Vert_{L^1(B)}^{p_0}} +1 \right]  \, \int_B | \d v |^{p_0} dg.
\ee Here
$v_f$ is the average of $v$ under weight $f$. i.e. $v_f =\int_B v
f dg/\int_B f dg.$ \proof
\end{lemma}
Under the basic assumptions, as shown in Section 2 of
\cite{TZq:1}, the following unweighted $L^2$ Poincar\'e inequality
is true:  for all smooth $v$ on $B$
\be \int_{B} |v-  v_B |^2 dg
\le C r^2 \int_{B} |\d v|^2 dg.
\ee Here $v_B$ is the average of $v$ on the ball $B$.
Since the manifold is volume
doubling, using the general result of \cite{KZ:1}, one can find a
positive number $p_0<2$ such that the unweighted $L^{p_0}$
Poincar\'e inequality holds: for all smooth $v$ on $B$ \be
\lab{lpPI} \int_{B} |v-  v_B |^{p_0} dg \le C r^{p_0} \int_{B} |\d
v|^{p_0} dg. \ee Here the constant $C$ may have changed. Let us
mention that the self improvement property of the Poincar\'e
inequalities was stated for the so called $1-q$ type, which means
the left hand side is in $L^1$ norm. However, as pointed out in
that paper, this property also holds for the above $2-2$ type
Poincar\'e inequality.

By the triangle inequality
\be
\lab{v-vf}
\int_B |v-v_f|^{p_0} dg \le 2^{p_0} \left( \int_B |v-v_B|^{p_0} dg +
 \int_B |v_B-v_f|^{p_0} dg \right).
\ee Observe that the second term inside the parentheses on the right hand side
satisfies
\be
\al
\int_B |v_B-v_f|^{p_0} dg &= |v_B-v_f|^{p_0}  |B|\\
&=\left| \frac{\int_B v_B f dg -\int_B v f dg}{\int_B f dg} \right|^{p_0} \, |B|
= |B| \left(\int_B f dg\right)^{-p_0} \left(\int_B (v-v_B) f dg \right)^{p_0}\\
&\le|B| \Vert f \Vert^{-p_0}_{L^1(B)} \Vert f \Vert^{p_0}_{L^\infty (B)}
\int_B |v-v_B|^{p_0}  dg \, |B \cap supp \, f|^{p_0-1}.
\eal
\ee Substituting this to the right hand side of (\ref{v-vf}) and using
({\ref{lpPI}}), we deduce
\be
\int_{B} |v-  v_f |^{p_0} dg
\le C \left[ |B| \, |B \cap supp \, f|^{p_0-1} \Vert f \Vert^{-p_0}_{L^1(B)}
\Vert f \Vert^{p_0}_{L^\infty (B)} + 1 \right]
 r^{p_0} \int_B |\d v|^{p_0}  dg,
\ee which proves the lemma.
\qed

The next lemma is a rerun of De Giorgi's method in our situation.
We will closely follow the presentation in \cite{Lie:1} Chapter
VI, Section 12. Since there are certain differences from the
Euclidean setting, due to the lack of local $L^1$ Poincar\'e
inequality, we will present a detailed proof.  First let us define
a De Giorgi class, which is not the most general one, but which is
sufficient for our setting.

\begin{definition}
\lab{deDeGi}
Let $D$ be a domain in $\M$. A function $v \in W^{1, 2}(D)$ is said to be in the
De Giorgi class in $D$ if the following inequalities hold. For all numbers $\sigma \in (0, 1)$,
$k$ and $B(p, r) \subset D$,
\be
\int_{B(p, \sigma r)} |\d (v-k)^+|^2 dg \le \frac{\beta_1}{(1-\sigma)^2 r^2}
\int_{B(p, r)} |(v-k)^+|^2 dg +  \beta_0 |D^+_k|,
\ee
\be
\int_{B(p, \sigma r)} |\d (v-k)^-|^2 dg \le \frac{\beta_1}{(1-\sigma)^2 r^2}
\int_{B(p, r)} |(v-k)^-|^2 dg +  \beta_0 |D^-_k|.
\ee Here $\beta_0$, $\beta_1$ are positive constants, $(v-k)^+=\max \{v-k, 0\}$,
$D^+_k=\{x \in B(p, r) \, | \, v-k > 0 \}$;  $(v-k)^-=\max \{-(v-k), 0\}$,
$D^-_k=\{x \in B(p, r) \, | \, v-k < 0 \}$.
\end{definition}

\begin{lemma}
\lab{leDeGi}
Suppose a function $v$ is in the De Giorgi class in $B(p, 2 r) \subset M$. Then there exist
a number $\a \in (0, 1)$ and a positive constant $C$, depending only on $\beta_0$, $\beta_1$
and the constants in the basic assumptions, such that, for $\rho \in (0, r)$,
\be
Osc_{B(p, \rho)} v \le C [(\rho/r)^\a Osc_{B(p, r)} v  + \beta_0 \rho].
\ee Here $Osc_D v$ is the oscillation of the function $v$ in the region $D$.
\proof
\end{lemma}
The proof is divided into three steps.

{\it Step 1.} We prove the following claim:
{\it Let $v \ge 0$ be in the De Giorgi class in $B(p, 2 r)$ and let $\tilde v = v + \beta_0 r$.
Suppose there exists  positive numbers $a<1$ and $K$ such that
\be
|\{ x \in B(p, r) \, | \, \tilde v <K\} | \le a |B(p, r)|.
\ee  Then, for any small number $\e \in (0, 1)$, there exists a number $\delta \in (0, 1)$
such that
\be
|\{ x \in B(p, r) \, | \, \tilde v < \delta K\} | \le \e |B(p, r)|.
\ee  Here $\delta$ depending only on $b$ and the constants in the statement of
 the lemma.} The claim says that if $\tilde v$ is larger than $K$ in a portion of the ball,
then the set where $\tilde v$ is much smaller than $K$ has a small measure.

To prove the claim, we introduce two functions on $B(p, r)$.
\be
w=w(x) =
\begin{cases}
0, \quad \text{if} \quad \tilde v(x) \ge a^i K,\\
a^i K - \t v(x) \quad \text{if} \quad a^i K \ge \t v(x) > a^{i+1} K,\\
a^i K -a^{i+1} K  \quad \text{if} \quad \t v(x) \le a^{i+1} K.
\end{cases}
\ee Here $i$ is a nonnegative integer.
\be
f=f(x) =
\begin{cases}   |\{x \in B(p, r) \, | \, w(x)=0 \}|^{-1} \quad \text{if} \quad
w(x)=0, \quad i.e. \quad \t v \ge a^i K,\\
0 \quad \text{if} \quad w(x) \neq 0, \quad i.e. \quad \t v < a^i K.
\end{cases}
\ee We also denote by $A_i$ the subset of $B(p, r)$ on which $\t v < a^i K$.
Thus $A_i^c$ is the set where $w$ vanishes, and $A_i$ is where $f$ vanishes.
Our assumption implies $|A_i^c| \ge (1-a) |B(p, r)|$ and also
$w f =0$. So the weighted average $w_f=0$, $\Vert f \Vert_{L^1(B(p, r))} =1$,
$\Vert f \Vert_{L^\infty(B(p, r))} = |A_i^c|^{-1}$.
From these we can apply Lemma \ref{lewfPI} to infer
\be
\al
\int_{B(p, r)} w^{p_0} dg
&\le C r^{p_0} \sup f^{p_0} \frac{|B| |B \cap supp \,
f|^{p_0-1}}{\Vert f \Vert_{L^1(B)}^{p_0}} \int_B | \d w |^{p_0} dg\\
&\le C r^{p_0}  \frac{|B| }{|A_i^c|} \int_{A_i-A_{i+1}} | \d \t v |^{p_0} dg,
\eal
\ee which shows
\be
[ a^i (1-a) K]^{p_0} |A_{i+1}| \le \int_{B(p, r)} w^{p_0} dg
\le C r^{p_0}  \int_{A_i-A_{i+1}} | \d \t v |^{p_0} dg
\ee By H\"older inequality, this implies
\be
\lab{ai1-aK}
[ a^i (1-a) K]^{p_0} |A_{i+1}| \le  C r^{p_0}
 \left( \int_{A_i} | \d \t v |^2 dg \right)^{p_0/2} \, |A_i - A_{i+1}|^{(2-p_0)/2}.
\ee From the assumption that $v$ is in De Giorgi class in $B(p, 2r)$, we know,
for all real numbers $k$,
\be
\int_{B(p,  r)} |\d (v-k)^-|^2 dg \le \frac{\beta_1}{r^2}
\int_{B(p, 2 r)} |(v-k)^-|^2 dg +  \beta_0 |D^-_k|.
\ee Here $D^-_k=\{x \in B(p, 2r) \, | \, v-k < 0 \}$. Hence
\be
\al
& \int_{B(p,  r)} |\d (\t v-k-\b_0 r)^-|^2 dg\\
&=\int_{B(p,  r)} |\d (v-k)^-|^2 dg
\le \frac{C }{r^2}
\int_{\{x \in B(p, 2r) \, | \, v(x)<k\}}[ |k - v|^2 + (\b_0 r)^2] dg\\
&\le \frac{C }{r^2} \int_{\{x \in B(p, 2r) \, | \, \t v(x)<k+\b_0
r \}}( |k+ \b_0 r - \t v|^2 + \t v^2) dg. \eal \ee In the last
step, we have used the fact that $\t v = v + \b_0 r \ge \b_0 r$.
Taking $k$ so that $k + \b_0 r = a^i K$, we then deduce \be
\int_{A_i} | \d \t v |^2 dg \le C (a^i K)^2 r^{n-2}. \ee Here we
have used the volume upper bound of geodesic balls. Substituting
this to the right hand side of (\ref{ai1-aK}), we find that \be [
a^i (1-a) K]^{p_0} |A_{i+1}| \le  C (a^i K)^{p_0} r^{n p_0/2} |A_i
- A_{i+1}|^{(2-p_0)/2}. \ee  Cancelling the term $(a^i K)^{p_0}$
and writing $\theta_i = r^{-n} |A_i|$, we infer \be
\theta^{2/(2-p_0)}_{i+1} \le C (1-a)^{-2p_0/(2-p_0)} (\theta_i -
\theta_{i+1}). \ee Notice that $\theta_i \ge \theta_{i+1}$. By
adding the above inequality from $i=0$ to a positive integer
$j-1$,  we conclude that \be j \theta^{2/(2-p_0)}_{j} \le C
(1-a)^{-2p_0/(2-p_0)} \theta_0 \le C, \ee which means \be |\{ x
\in B(p, r) \, | \, \t v(x) < a^j K \} | \le C j^{-1} |B(p, r)|.
\ee The claim is proven by choosing $j$ large enough so that $C
j^{-1} \le \e$ and write $\delta =a^j$.
\medskip

{\it Step 2.} We prove the following assertion: {\it there exists a number $\e \in (0, 1)$
such that if,
\be
|\{ x \in B(p, r) \, | \, \t v(x) < K \} | \le \e |B(p, r)|
\ee for some $K>0$, then $\t v \ge K/2 - C_0 r$ on $B(p, r/2)$.
Here $C_0$ is a constant depending only on the basic assumptions for $\M$.}

Following the proof in the Euclidean case verbatim, since $v$ is in the De Giorgi class,
we know that the following mean value
inequality is true:
\be
\sup_{B(p, r/2)} [(K -\t v)^+]^2 \le \frac{C}{r^n} \int_{B(p, r/2)}
[(K- \t v)^+]^2 dg + C (\b_0 r)^2.
\ee The proof just uses the $L^2$ Sobolev inequality and volume non-collapsing
property of geodesic balls. We refer the reader to p143 of \cite{Lie:1}
for a detailed proof. Note that the proof there is for the more general parabolic
case. The above mean value inequality shows
\be
\sup_{B(p, r/2)} (K -\t v)^2 \le \frac{C}{r^n} |\{ x \in B(p, r) \, | \, \t v(x) < K \} |
K^2
+  C (\b_0 r)^2 \le C \e K^2 + C (\b_0 r)^2.
\ee If $C \e<1/8$ then for all $x \in B(p, r/2)$, the preceding inequality implies
\be
\t v(x) \ge K/2 - C_0 r,
\ee which proves the assertion.
\medskip

{\it Step 3.} We will finish the proof of the lemma.

Given a number $\rho \in (0, 2r]$, we write $M(\rho) = \sup_{x \in B(p, \rho)} v(x)$,
$m(\rho) = \inf_{x \in B(p, \rho)} v(x)$ and $J(\rho) = M(\rho)-m(\rho)$ i.e. the
oscillation of $v$ in $B(p, \rho)$.
Consider the function  $h=h(x) = v(x) - m(2 r) $. Then $h$ is a nonnegative function in the
De Giorgi class in $B(p, 2 r)$.
Note that
\be
\{ x \in B(p, r) \, | h(x) \le J(2r)/2 \} \cup \{ x \in B(p, r) \, | h(x) > J(2r)/2 \}
=B(p, r).
\ee Thus we can assume, without loss of generality, that
\be
| \{ x \in B(p, r) \, | h(x) \le J(2 r)/2 \}| \le \frac{1}{2} |B(p, r)|,
\ee since we can consider $J(2 r) -h=M( 2r) -v$ otherwise.
Write $\t h = h + \b_0 r$. Then
\be
\lab{j2r+b0r}
| \{ x \in B(p, r) \, | \t h(x) \le J(2 r)/2 + \b_0 r \}| \le \frac{1}{2} |B(p, r)|,
\ee  Let $A>8 \b_0 $ be a large positive number to be specified later.
If $J(2 r) \le A r$, then we stop and rerun the above process on the ball $B(p, r)$
instead. So we assume $J(2r) > A r$.  In this case, (\ref{j2r+b0r}) implies
\be
| \{ x \in B(p, r) \, | \t h(x) \le 3 J(2 r)/4 \}| \le \frac{1}{2} |B(p, r)|.
\ee By Step 1, for any $\e>0$, there is $\delta>0$, such that
\be
| \{ x \in B(p, r) \, | \t h(x) \le \delta 3J(2 r)/4 \}| \le \e |B(p, r)|.
\ee Choose $\e$ to be the number in the assertion in Step 2 and $K=\delta 3 J(2r)/4$.
Then, for $x \in B(p, r/2)$, we have
\be
\t h(x) \ge  \delta 3J(2 r)/8 - C_0 r \ge \delta J(2 r)/4.
\ee In the last inequality we used $J(2r) \ge A r$ and choose $A$ sufficiently large.
Therefore
\be
v(x)  \ge \delta J(2r)/4 -\b_0 r + m(2r)  \ge \delta J(2r)/8 + m(2r).
\ee i.e.
\be
m(r/2) \ge \inf_{x \in B(p, r/2)} v(x)  \ge \delta J(2r)/8 + m(2r).
\ee This implies
\be
J(r/2) = M(r/2) - m(r/2) \le M(2r) - m(2r) - \frac{\delta}{8} J(2r)
= (1-\frac{\delta}{8}) J(2r).
\ee To summarize, we have proven that either $J( 2r) \le A r$ or
$J(r/2) \le (1-\frac{\delta}{8}) J(2r)$. Repeating this on the balls $B(p, r)$, $B(p, r/2)$, ..., we
have proven the lemma is true.
\qed

Since we are concerned with only local properties when applying the above lemma,
we will take the radius $r \le 1$ in the rest of the section.

\begin{lemma}
\lab{lec1aharm} Suppose the ball $B(p, 2r)$ is contained in a
harmonic coordinate $\{x^1, ..., x^n\}$. We assume that the metric satisfies, in the ball,
the following $C^{1}$ bound.
a). $  e^{-\theta} I \le (g_{pq}) \le e^\theta  I , \, I = (\delta_{pq})$;
b). $  \sup_{p, q} ( r \, \Vert g_{pq} \Vert_{C^{1}}  ) \le e^{\theta}$.
Here $\theta>0$.

Let $h$ be a harmonic
function in $B(p, r)$ and $dh = h_i dx^i$. Then exist positive
constants $\a_0 \in (0, 1)$ and $C$, which depend only on the
parameters of the basic assumptions and $\theta$, such that \be \lab{c1ah}
\Vert  h_i \Vert_{C^{\a_0}(B(p, r))} \le C \left( \frac{1}{r^{\a_0}}
\Vert \d h \Vert_{L^\infty(B(p, 2r))} + 1 \right).
\ee
\proof
\end{lemma}

We work in the real coordinate system. Since $\Delta h =0$, we know that $d h$
is a harmonic one form, i.e. $d d^* d h=0$. In the harmonic system, we write
$d h = h_i dx^i$. Since $x^i$ is a harmonic function, we know, for any constants
$k_i$, that
\be
(d d^* + d^* d) (  (h_i - k_i) dx^i) =0.
\ee From this, the Weitzenboch formula implies
\be
\lab{eqhiki}
\Delta (h_i-k_i) - R_{ij} (h_j-k_j) =0
\ee where $\Delta$ is the rough Laplacian on $1$ forms.
Denoting  $\eta=\eta_i dx^i$ for the one form $(h_i-k_i) dx^i$, we will write down
 equation (\ref{eqhiki}) for $\eta_i$ in the local system.

Let $\d_i$ be the covariant derivative in the $x^i$ direction, and
$\p_i$ be the partial derivative. Then \be \d_i \eta = (\p_i
\eta_k - \eta_l \Gamma^l_{ik}) dx^k. \ee Here and later in this
paragraph, with a slight abuse of notation, we use $k$ as an index
rather than the free constant in the De Giorgi argument. But it
should be clear from the context what $k$ means.  Further more \be
\d_j (\d_i \eta) = \left[ \p_j  (\p_i \eta_k - \eta_l
\Gamma^l_{ik}) - (\p_i \eta_m -\eta_l \Gamma^l_{jm}) \Gamma^m_{jk}
\right] dx^k . \ee Recall that the local formula for the rough
Laplacian of $\eta$ is \be g^{ij} \d_j (\d_i \eta) +
\frac{1}{\sqrt{det g}} \p_i (\sqrt{det g}\, g^{ij})\, \d_j \eta.
\ee See Section 10.1 of \cite{Ni:1} e.g. In harmonic coordinates,
it then takes the form \be \Delta \eta = g^{ij} \d_j (\d_i \eta)
=g^{ij} \left[ \p_j  (\p_i \eta_k - \eta_l \Gamma^l_{ik}) - (\p_i
\eta_m -\eta_l \Gamma^l_{jm}) \Gamma^m_{jk} \right] dx^k. \ee
Substituting this identity into (\ref{eqhiki}), we find that \be
g^{ij} \p_i \p_j \eta_k - g^{ij} \p_j (\eta_l \Gamma^l_{ik})
 - g^{ij} \p_i \eta_m \Gamma^m_{jk} + g^{ij} \eta_l \Gamma^l_{im}
\Gamma^m_{jk} =  g^{lm} R_{lk} \, \eta_m. \ee Taking $\eta_k$ as a
scalar function and $\Delta$ be the scalar Laplacian, then we
have, in the harmonic coordinates \be \lab{scalaDD} \Delta \eta_k
- g^{ij} \p_j (\eta_l \Gamma^l_{ik})
 - g^{ij} \p_i \eta_m \Gamma^m_{jk} + g^{ij} \eta_l \Gamma^l_{im} \Gamma^m_{jk}
-  g^{lm} R_{lk} \, \eta_m =0.
\ee

Let $\phi$ be a Lipschitz cut off function supported in $B(p, 2r)$
such that $\Vert \d \phi \Vert_\infty \le C/r$. Next we fix the
index $i$ and take $k_i $ to be a free constant $k$, and $k_j=0$
for $j \neq i$.  Denote by $(h_i-k)^+$ the positive part of
$h_i-k$. Using $(h_i -k)^+ \phi^2$ as a test function in the
equation (\ref{scalaDD}), after switching the indices suitably and
doing integration by parts, we find that
\be \al
\int &|\d [(h_i-k)^+ \phi]|^2 dg \\
&\le C [\Vert \d \phi \Vert^2_\infty+1]
\int_{supp \, \phi} [(h_i-k)^+]^2 +  C \Vert \eta \Vert_{L^\infty}  \,
\Vert \Gamma^m_{jl} \Vert_{L^\infty} \int |\d [(h_i-k)^+ \phi]| dg\\
&\qquad + C \Vert \Gamma^m_{jl} \Vert_{L^\infty}
\int  (h_i-k)^+ \phi \,  | Hess \, h| \phi dg+ \Vert \eta \Vert_{L^\infty}
C \Vert \Gamma^m_{jl} \Vert^2_{L^\infty} \, \int  [(h_i-k)^+ \phi] dg \\
& \qquad  + C \Vert \eta \Vert_{L^\infty} \int |Ric|  (h_i - k)^+
\phi^2 dg. \eal
\ee
Here $\d$ stands for the gradient of a scalar
function, $\Vert \Gamma^m_{jl} \Vert_{L^\infty}$ is the $L^\infty$
norm of Christoffel symbols. We have also used the fact
\be
| g^{ij} \p_i \eta_m | \le C | Hess \, h | + C e^{\theta} | \d h|.
\ee
By the assumed $C^{1}$ bound for the
metric, for any $\e>0$, we deduce
\be
\al
\int &|\d [(h_i-k)^+
\phi]|^2 dg \\
&\le C (\Vert \d \phi \Vert^2_\infty +1) \int_{supp \,
\phi} [(h_i-k)^+]^2  + \e \int [|Ric|^2 + | Hess \, h
|^2]  [ (h_i - k)^+ \phi]^2 dg \\
&\qquad+ C \e^{-1} \int_{D^+_k} \Sigma_{j
\neq i} (h_j-k_j)^2 dg. \eal \ee Here $D^+_k = \{x | \, h_i - k
\ge 0 \}$. Using  the embedding Lemma \ref{leimbed}, we see that
\be \al
&\int |\d [(h_i-k)^+ \phi]|^2 dg \\
&\le C (\Vert \d \phi \Vert^2_\infty +1)
\int_{supp \, \phi} [(h_i-k)^+]^2 +
C \e K(|Ric|^2+| Hess \, h |^2) \int |\d  [ (h_i - k)^+ \phi]|^2 dg \\
&\qquad+ C \e K(|Ric|^2+| Hess \, h |^2) r^{-2} \int   [ (h_i - k)^+ \phi]^2 dg +
C \e^{-1} \Vert \d h \Vert^2_\infty |D^+_k|.
\eal
\ee Here $K(|Ric|^2+| Hess \, h |^2)$ is defined in
(\ref{Khess}), which is a bounded quantity by Lemma \ref{lehess2}.  Although
 we only proved $K(|Hess \, u|^2)$ is bounded, the proof for the boundedness of
$K(|Hess \, h|^2)$ is the same, and actually simpler.  Here $u$ is
the Ricci potential.

Choosing
$\e = [ 2 C K(|Ric|^2+| Hess \, h |^2)]^{-1}$, we reach
\be
\lab{hik+}
\al
\int |\d [(h_i-k)^+ \phi]|^2 dg
&\le C \left(\Vert \d \phi \Vert^2_\infty +\frac{K^2(|Ric|^2+| Hess \, h |^2)}{r^2}+1 \right)
\int_{supp \, \phi} [(h_i-k)^+]^2 \\
&\qquad+
C K^2(|Ric|^2 +| Hess \, h |^2)  \Vert \d h \Vert^2_\infty |D^+_k|.
\eal
\ee In the same manner, we deduce, for $(h_i-k)^- = - \min \{0, h_i-k \}$, that
\be
\lab{hik-}
\al
\int |\d [(h_i-k)^- \phi]|^2 dg &\le C \left(\Vert \d \phi \Vert^2_\infty
+\frac{K^2(|Ric|^2+| Hess \, h |^2)}{r^2}+1 \right) \int_{supp \, \phi}
[(h_i-k)^-]^2 \\
&\qquad+ C K^2(|Ric|^2+| Hess \, h |^2)  \Vert \d h \Vert^2_\infty |D^-_k|.
\eal \ee Here $D^-_k = \{x | \, h_i - k \le 0 \}$. Inequalities
(\ref{hik+}) and (\ref{hik-}) mean that $h_i$ is in the
De Giorgi class.
By Lemma \ref{leDeGi}
we know that $h_i$ is
$C^{\a_0}$ for some $\a_0 \in (0, 1)$ and (\ref{c1ah}) holds.  \qed

\begin{lemma}
\lab{lec1apot} Suppose the ball $B(p, 2r)$ is contained in a
harmonic coordinate $\{z^1, ..., z^m \}$. Assume the same $C^{1}$ bound
for the metric as in the previous lemma.
Let $u$ be the Ricci
potential and $u_i$ be the component of $du$ with respect to
$dz_i$. Then exist positive constants $\a_0 \in (0, 1)$ and $C$,
which depend only on the parameters of the basic assumptions and the $C^1$
bound for the metric, such
that
\be \Vert  u_i \Vert_{C^{\a_0}(B(p, r))} \le C \left(
\frac{1}{r^{\a_0}}  \Vert \d u \Vert_{L^\infty(B(p, 2r))} + 1 \right).
\ee \proof
\end{lemma}
The proof is similar to the previous lemma. One difference is that in complex coordinates, the derivative of $u$ is complex valued. In order to apply the De Giorgi
method we need to treat the real and imaginary parts separately.  Another difference
is an extra term involving the scalar curvature.
Let $\{z^1, ..., z^m \}$ be a complex harmonic coordinate. Recall that
$\Delta u = m-R$. Write $du = u_i dz^i$ and let $k_i$ be real constants. Then
\be
\Delta (u_i-k_i)= \frac{1}{2} R_{i \g j} (u_j - k_j) - R_i, \quad
\Delta  (u_{\g i}-k_i)= \frac{1}{2} R_{\g i  j} (u_{\g j} - k_j) - R_{\g i},
\ee Adding these two equations together, we see that
\be
\lab{eqreuik}
\Delta (Re \, u_i - k_i) = \frac{1}{4} R_{i \g j} (u_j - k_j)+
\frac{1}{4} R_{\g i  j} (u_{\g j} - k_j) -  Re \, R_i
\ee As in the previous lemma, we fix an index $i$ and take $k_i$ to be a real
constant $k$, and $k_j=0$ for all $j \neq i$. Pick a Lipschitz cut off function
$\phi$ supported in $B(p, 2r)$ such that $|\d \phi| \le C/r$. Using $(Re \, u_i - k_i)^+
\phi^2$ as a test function in (\ref{eqreuik}), after carrying out the same local calculation
as the previous lemma, we deduce
\be
\al
&\int |\d [(Re \, u_i - k)^+ \phi] |^2 dg \\
&\le C (\Vert \d \phi \Vert^2_\infty+1)
\int_{D_k}  [(Re \, u_i - k)^+ \phi] |^2 dg + \int (|R_{i \g j}| + |Hess \, u|) |u_j-k_j|
(Re \, u_i - k)^+ \phi^2 dg \\
&\qquad+ \left| \int Re \, R_i  (Re \, u_i - k)^+ \phi^2 dg
\right|.
\eal
\ee Doing integration by parts on the last term and Using Cauchy-Schwarz inequality, we arrive at
\be
\al
&\int |\d [(Re \, u_i - k)^+ \phi] |^2 dg \\
&\le C ( \Vert \d \phi \Vert^2_\infty + 1)
\int_{D_k}  [(Re \, u_i - k)^+ \phi] |^2 dg + \int |Hess \, u|^2
[(Re \, u_i - k)^+ \phi]^2 dg \\
&\qquad+ \int_{D^+_k} \Vert \d u \Vert^2_\infty \phi^2 dg+ 2 \int | R|  |\d [(Re \, u_i - k)^+ \phi^2]| dg.
\eal
\ee Here $D^+_k = \{ x \in B(p, 2 r) \, | \, Re \, u_i(x) - k \ge 0 \}$.
By the embedding
Lemma \ref{leimbed}, we can turn the above inequality into
\be
\al
&\int |\d [(Re \, u_i - k)^+ \phi] |^2 dg \\
&\le C ( \Vert \d \phi \Vert^2_\infty +\Vert \d \phi \Vert^2_\infty  \Vert R \Vert^2_\infty+ 1)
\int_{D_k}  [(Re \, u_i - k)^+ \phi] |^2 dg\\
&\qquad \qquad +  C \left[ K^2(|Hess \, u|^2) (\Vert \d u \Vert^2_\infty +1) + \Vert R \Vert^2_\infty \right] \, |D^+_k|.
\eal
\ee  Likewise, we also have
\be
\al
&\int |\d [(Re \, u_i - k)^- \phi] |^2 dg \\
&\le C ( \Vert \d \phi \Vert^2_\infty +\Vert \d \phi \Vert^2_\infty  \Vert R \Vert^2_\infty+ 1)
\int_{D_k}  [(Re \, u_i - k)^- \phi] |^2 dg\\
&\qquad \qquad +  C \left[ K^2(|Hess \, u|^2) (\Vert \d u \Vert^2_\infty +1) + \Vert R \Vert^2_\infty \right] \, |D^-_k|.
\eal
\ee Therefore, $Re \, u_i$ is in the De Giorgi class. Similarly, the same is true for
$Im \, u_i$.
 By Lemma
\ref{leDeGi} again, we now
know that $u_i$ is H\"older continuous, proving the lemma.
\qed

\section{lower bound of harmonic radius}

In the previous section we proved that the gradient of a harmonic function is H\"older
continuous within a harmonic chart. In this section we will show that the harmonic chart
can not be too small under certain conditions.

\begin{definition} (maximal harmonic radius)
\lab{defharmrad}
Given  numbers $\theta>0$ and  $\a \in (0, 1)$,  and a point $x \in \M$ with metric $g$.
 We define
$r^{\theta, \a}_g(x)$ to be the maximum radius $l$ such that there exists,
in the ball $B(x, l)$,  a $C^{1, \a}$ harmonic
coordinate  $X=(x^1, .., x^n) :
B(x, l) \to {\bf R}^n$,  which satisfies the following properties.

\be \al \lab{th-alcond}
 & a). \quad  e^{-\theta} I \le (g_{pq})
\le e^\theta  I , \, I = (\delta_{pq}); \\
& b). \quad  \sup_{p, q}  l \Vert \d g_{pq} \Vert_{C^0} \le
e^\theta \\
& c). \quad \sup_{p, q} l^{1 + \a}
 \Vert \d g_{pq} \Vert_{C^{\a}}  \le e^{\theta}.
\eal \ee
\end{definition}

\begin{lemma}
\lab{leharmrad}
Let $({\M}, g)$ be a manifold satisfying the basic assumptions, and $x$ be a point
in $\M$.
Let $\theta>0$ be any fixed number and
$\a \in (0, \a_0)$ where $\a_0$ is the H\"older parameter in Lemmas
\ref{lec1aharm} and
\ref{lec1apot}.

Suppose either one of the following conditions holds.

1.  the injectivity radius at  each point of $B(x, r)$ is greater than  a number $i_0>0$.

2. the geodesic ball $B(x, r)$ is almost Euclidean, i.e.  for a sufficiently
small positive number $\eta$, any $y \in B(x, r)$ and
$0<\rho <dist(y, \partial B(x, r))$, it holds
\be
|B(y, \rho)| \ge (1-\eta) \omega_n \rho^n
\ee where $\omega_n$ is the volume of the unit ball in ${\bf R}^n$.

3. for a sufficiently small positive number $\eta>0$ and a given number $\tau_0>0$,
 the infimum of the  $W$ entropy on $B(x,  r)$
satisfies
\be
\sup_{\tau \in (0, \tau_0]} \mu(g, \tau, B(x, r)) \ge -\eta.
\ee

 Then there exists a number $\delta \in (0, 1)$ such that
\be
r^{\theta, \a}_g (y) \ge \delta dist (y, \partial B(x, r))
\ee for all $y \in B(x, r)$.

In particular the harmonic radius at $x$ satisfies
\be
r^{\theta, \a}_g (x) \ge \delta r.
\ee
\end{lemma}

\proof

We will use the blow up method to prove the lemma.
The strategy follows that in the paper \cite{An:1} Main Lemma 2.2, where
Anderson proved the same conclusion under the assumption that the Ricci curvature
is bounded and the injectivity radius lower bound.  In the current situation, we do not know whether the Ricci curvature is bounded. To overcome this difficulty, we
will use the regularity result proved in the previous lemmas. Another difference
occurs since we assume condition 3 instead of the lower bound on the injectivity
radius.

We divide the rest of the proof into three cases.

{\it Case 1.} Suppose Condition 1 holds.

The proof of this case consists of most of the work for the lemma, on which
the proof of the other 2 cases are built upon.

Assuming the conclusion of the lemma is false, then there exists a sequence of pointed manifolds
$({\M}_i, x_i, g_i)$ which satisfies the basic assumptions, but, for some
$p_i \in B(x_i, r, g_i) \subset {\bf M}_i$, we have
\be
\al
\frac{r^{\theta, \a}_{g_i}(p_i) }{dist( p_i, \, \p B(x_i, r, g_i))}
&=  \inf \{ \frac{r^{\theta, \a}_{g_i}(y) }{dist( y, \, \p B(x_i, r, g_i))}  \, | \,
 y \in B(x_i, r, g_i)
 \} \\
&=  \delta_i \to 0.
\eal
\ee Notice that $p_i$ are the worst points in the sense that for every other point
$y \in B(x_i, r, g_i)$,
we have
\be
\frac{r^{\theta, \a}_{g_i}(y)}{dist( y, \, \p B(x_i, r, g_i))}
 \ge \frac{r^{\theta, \a}_{g_i}(p_i) }{dist( p_i, \, \p B(x_i, r, g_i))}.
\ee

For simplicity, we write
\be
\lambda_i = dist( p_i, \, \p B(x_i, r, g_i)).
\ee
From the Definition \ref{defharmrad}, the following statement holds: if $\rho>
\delta_i \lambda_i$, then
there are no harmonic coordinate systems
satisfying bounds (\ref{th-alcond}) on balls $ B(p_i, \rho, g_i)$.

Now, consider the scaled metric
\be
h_i = (\delta_i \lambda_i)^{-2} g_i.
\ee
Then by the above discussion, for every point
$y \in B(p_i,  \delta^{-1}_i, h_i)$, the ball $B(y, 1, h_i)$ is contained in a harmonic
coordinate such that  on this same ball and under $h_i$, we have
\be
\lab{hith-alcond}
  e^{-\theta} I \le ((h_i)_{pq}) \le e^\theta  I  \quad \text{and}
 \quad  \sup_{p, q} ( \Vert (h_i)_{pq} \Vert_{C^{1, \a}}  ) \le e^{\theta}.
\ee However, for any number $\rho>1$,
\be
\lab{nolargerrad}
\text{there is no harmonic coordinate system
containing} \quad B(p_i, \rho, h_i)
\ee which satisfies (\ref{hith-alcond}) on $B(p_i, \rho, h_i)$ under the metric $h_i$.

By Lemma 2.1 in \cite{An:1}, there is a subsequence of the triple
$\{ B(p_i,  \delta^{-1}_i, h_i), p_i, h_i \}$, still denoted by the same notation,
which converges in $C^{1, \a'}_{loc}$ topology to a complete, pointed manifold $(N, z, h)$. Here $\a'$ is any positive
number strictly less than $\a$. Also $h$ is a $C^{1, \a}$ metric.  Since the original
manifold has uniformly $C^1$ Ricci potential $u_i$, it is easy to see that the limit manifold is Ricci
flat. The reason is that the components of Ricci curvature  $R_{p \g q}$ does not change after scaling so, under the
metric $h_i$, it becomes
\be
R_{p \g q} =- \p_p \p_{\g q} u_i + (\delta_i \lambda_i)^2 (h_i)_{p \g q}.
\ee But
\be
| \d_{h_i} u_i| = \delta_i \lambda_i | \d_{g_i} u_i| \to 0
\ee as $i \to \infty$. Therefore the limiting manifold is Ricci flat since its Ricci potential
 is a constant.

 Condition 1 then implies that the injectivity radius at $z$ is infinity under the
metric $h$. The Cheeger-Gromoll splitting theorem tells us that $(N, z, h)$ is in fact
the standard Euclidean space.

 We claim that, by choosing a subsequence suitably,
we can ensure that the above convergence is actually in the original $C^{1, \a}_{loc}$
topology.
Taking the claim for granted, then, for the fixed $\theta$ and any fixed $\rho>1$, we know that (\ref{hith-alcond}) holds in $B(p_i, \rho, h_i)$
when $i$ is sufficiently large. This is a contradiction to (\ref{nolargerrad}), which would prove the lemma.

 Since $\a < \a_0$, bounded sets in $C^{1, \a_0}_{loc}$
topology are compact in $C^{1, \a}_{loc}$ topology.
Hence the claim  holds if we can prove that, for any ball
$B(w, 1, h_i) \subset  B(p_i, 2 \delta^{-1}_i, h_i)$, the metric $h_i$ satisfies, for a number
$\a_1 \in (\a, \a_0)$, that
\be
\lab{hic1a0}
  \sup_{p, q} [ \Vert (h_i)_{pq} \Vert_{C^{1, \a_1}(B(w, 1/2, h_i))}  ] \le  C
\ee for a uniform constant $C=C(\a_1)$.  So, in order to verify the lemma,  we just need to
 prove (\ref{hic1a0}).

For the rest of the proof, we will suppress the index $i$ in the metric $h_i$.
Since $B(w, 1, h)$ is contained in a harmonic coordinate, say, $\{x^1, ..., x^n \}$, we know that the
components of $h$ satisfies
\be
\lab{eqhrs}
h^{jk} \frac{\p^2 h_{rs}}{\p x^j \p x^k} + Q(h_{pq}, \frac{\p h_{pq}}{\p x^m}) =
(Ric_h)_{rs},
\ee where $Q$ is a quadratic term involving th first covariant derivative of $h$;
$Ric_h$ is the Ricci curvature with respect to $h$. For simplicity we write the above
equation as
\be
\lab{eqhrs2}
\Delta  h_{rs}+ Q(h, \d h) =
R_{rs}.
\ee  Let $\phi$ be a Lipschitz cut-off function in $B(w, 1, h)$ such that $\phi =1$
in $B(w, 2/3, h)$ and $|\d \phi| \le 4$. Let $\Gamma$ be the Green's function on the whole manifold $(M, h)$. Using Green's formula, it is easy to see that
\be
\lab{eqdh}
\al
\d_{x_j} ( \phi h_{rs}(x)) &= \int \d_{x_j} \Gamma (x, y) Q(\d h) \phi(y) dg(y)
- 2 \int \d_{x_j} \Gamma(x, y) \d \phi \d h_{rs}(y) dg(y)\\
&\qquad - \int \d_{x_j} \Gamma (x, y)  h_{rs} \Delta \phi dg(y) -
\int \d_{x_j} \Gamma (x, y)  R_{rs}  \phi(y) dg(y)\\
&\equiv I_1 + I_2 + I_3 + I_4.
\eal
\ee For $x \in B(y, 1/2)$ we will prove that $I_1, ..., I_4$ are uniformly bounded in
$C^{\a_0}$ norm.  Here and later in the proof, when mentioning balls, distance etc we will suppress the reference to
the underlying metric $h$, unless stated otherwise.

By Lemma \ref{lec1aharm}, we know that for $z_1, z_2 \in B(w, 1/2)$ and $y \in B(w, 1)$, the following inequalities hold:  if $d(z_1, y) \le d(z_2, y)$ then
\be
\lab{dGz1y}
|\d_{x_j} \Gamma(z_1, y) - \d_{x_j} \Gamma(z_2, y)|
\le C [d(z_1, z_2)/d(z_1, y)]^{\a_0} \, d(z_1, y)^{-(n-1)};
\ee if $d(z_1, y) \le d(z_2, y)$, then
\be
\lab{dGz2y}
|\d_{x_j} \Gamma(z_1, y) - \d_{x_j} \Gamma(z_2, y)|
\le C [d(z_1, z_2)/d(z_2, y)]^{\a_0} \, d(z_2, y)^{-(n-1)},
\ee Here is the proof.  Assume $d(z_1, y) \le d(z_2, y)$. If also $d(z_1, y)
\ge 2 d(z_1, z_2)$, we just apply Lemma \ref{lec1aharm} on the ball $B(z_1, 2 d(z_1, y)/3)$, which says
\be
\al
&|\d_{x_j} \Gamma(z_1, y) - \d_{x_j} \Gamma(z_2, y)|\\
&\qquad \le C [d(z_1, z_2)/d(z_1, y)]^{\a_0} \sup_{z \in B(z_1, 2 d(z_1, y)/3)}
|\d_z \Gamma(z, y)| + C d(z_1, z_2)^{\a_0}.
\eal
\ee This implies (\ref{dGz1y}) by the gradient bound
on the Green's function.  On the other hand, if $d(z_1, y)
\le 2 d(z_1, z_2)$, then by this gradient bound  again, we
have
\be
|\d_{x_j} \Gamma(z_1, y)| \le C d(z_1, y)^{-(n-1)} \le
C [d(z_1, z_2)/d(z_1, y)]^{\a_0} \, d(z_1, y)^{-(n-1)},
\ee and, since $d(z_1, y) \le d(z_2, y)$ by assumption,
\be
|\d_{x_j} \Gamma(z_2, y)| \le C d(z_2, y)^{-(n-1)} \le C d(z_1, y)^{-(n-1)}
\le
C [d(z_1, z_2)/d(z_1, y)]^{\a_0} \, d(z_1, y)^{-(n-1)}.
\ee These two inequalities also imply (\ref{dGz1y}).
Similarly we can prove (\ref{dGz2y}).
 We remark that the
controlling parameters for the Ricci curvature etc do not become
worse when the metric is magnified. Recall that the metric $h$,
which is $h_i$ after restoring the index, is an magnification of
the original metric. Therefore we have uniform constants through
out. In the above we have also used assumption \ref{th-alcond} (a)
which allows us to bound, in the same time, the H\"older norm of
the gradient of $\Gamma$ and the components of $d \Gamma$ in the
coordinate.

We observe, due to (\ref{hith-alcond}), that $Q(h, \d h)$ is a bounded function on
$B(m, 1)$.
Hence, from (\ref{dGz1y}) and (\ref{dGz2y}), we see that
\be
\lab{i1ca}
\Vert I_1 \Vert_{C^{\a_0}(B(m, 1/2))} \le C \sup_{z \in B(m, 1/2)} \int_{B(m, 1)} d(z, y)^{-n+1-\a_0} dg(y) \le C,
\ee where we have used the volume non-inflating property.

Next  comes $I_2$.  Since $|\d \phi| \le 4$ and $|\d h_{rs}|$ is bounded by
(\ref{hith-alcond}) again, we know,  from (\ref{dGz1y}) and (\ref{dGz2y}) that
\be
\lab{i2ca}
\Vert I_2 \Vert_{C^{\a_0}(B(m, 1/2))} \le C \sup_{z \in B(m, 1/2)} \int_{B(m, 1)} d(z, y)^{-n+1-\a_0} dg(y) \le C.
\ee

Now we work on $I_3$ which should be written as
\be
\lab{i3=}
\al
I_3 &= \int h_{kl}(y) \d_{x_j} \d_{y_k} \Gamma(x, y) \d_{y_l} \phi(y) h_{rs}(y)
dg(y) \\
&\qquad + \int  \d_{x_j}  \Gamma(x, y) h_{kl}(y) \d_{y_l} \phi(y)
\d{y_k} h_{rs}(y) dg(y) \equiv I_{31}+I_{32}. \eal \ee  For fixed
$y$, as a the scalar function of $x$, we know that
\be \lab{dGy}
\d_{y_k} \Gamma(x, y) \equiv h(\d \Gamma(x, y), \frac{\p}{\p y_k})
\ee is a harmonic function of $x \neq y$. This can be seen by
picking a curve in $c=c(s)$ on the manifold, whose tangent at
$s=0$ is $\frac{\p}{\p y_k}$ and differentiating with respect to
$s$ on the identity \be \Delta_x \Gamma(x, c(s)) =0. \ee

For $y \in B(m, 2/3)^c$ and $z_1, z_2 \in B(m, 1/2)$, we can apply
Lemma \ref{lec1aharm} on $\d_{y_k} \Gamma(x, y)$ to deduce \be
\lab{ddGam} \al
&|\d_{x_j} \d_{y_k} \Gamma(z_1, y)- \d_{x_j} \d_{y_k} \Gamma(z_2, y)|\\
&
\le C \max \{ [d(z_1, z_2)/d(z_1, y)]^{\a_0} d(z_1, y)^{-n},
[d(z_1, z_2)/d(z_2, y)]^{\a_0} d(z_2, y)^{-n}
 \}\\
&\le C d(z_1, z_2)^{-\a_0}.
\eal
\ee Here we also used the bound
\be
|\d_{x_j} \d_{y_k} \Gamma(x, y)| \le \frac{C}{d(x, y)^n}
\ee which can be proven by using the gradient bound for the Green's function
twice on suitably chosen balls. We mention that $\d_{x_j} \d_{y_k} \Gamma(x, y)$
means the component in the $x_j$ direction of the gradient of the scalar function
$\d_{y_k} \Gamma(x, y)$.

 Thus $\Vert I_{31} \Vert_{C^{\a_0}(B(m, 1/2)} \le C$.
Similar to the case of $I_2$, we also have $\Vert I_{32} \Vert_{C^{\a_0}(B(m, 1/2)} \le C$. Therefore
\be
\lab{i3ca}
\Vert I_3 \Vert_{C^{\a_0}(B(m, 1/2))} \le C.
\ee

Finally we come to $I_4$.  Recall that $R_{rs} = - \p_r \p_s u +
h_{rs}$ where $r$ and $s$ implicitly mean unbarred and barred
index. Thus we can write, after integration by parts,
 \be
\lab{i4=} \al
I_4&= - \int \d_{x_j} \d_r \Gamma(x, y) (\p_s u(y) - \p_s u(q)) \phi(y) dg(y)\\
&\qquad -\int \d_{x_j} \Gamma(x, y)  (\p_s u(y) - \p_s u(q)) \d_r
\phi(y) dg(y)
 -\int \d_{x_j} \Gamma(x, y) h_{rs} \phi(y) dg(y)\\
 &\qquad -\int \d_{x_j} \Gamma(x, y)  (\p_s u(y) - \p_s u(q)) \phi
 \,
 \p_r (\ln \sqrt{det g}) dg(y)\\
&\equiv I_{41} + I_{42} + I_{43}+I_{44}, \eal \ee where $q$ is a
point to be chosen later. Also the scalar function $\d_r \Gamma$
is defined in (\ref{dGy}) and $\d_r \phi$ is similarly defined.

Recall that in  the ball $B(m, 1)$,  $|\d u|$ is uniformly
bounded, and by (\ref{hith-alcond}),
 \be
 \lab{thetajiehpq}
  e^{-\theta} I \le (h_{pq}) \le e^\theta  I  \quad \text{and}
 \quad  \sup_{p, q} ( \Vert h_{pq} \Vert_{C^{1, \a}}  ) \le e^{\theta}.
\ee We can now follow the case for $I_2$ to show that \be
\lab{i423ca} \Vert I_{42} \Vert_{C^{\a_0}(B(m, 1/2))} + \Vert
I_{43} \Vert_{C^{\a_0}(B(m, 1/2))} + \Vert I_{44}
\Vert_{C^{\a_0}(B(m, 1/2))}\le C. \ee So we are left with treating
$I_{41}$.

Pick $z_1, z_2 \in B(m, 1/2)$ which is divided into two regions
\be B(m, 1/2) = D_1 \cup D_2 \equiv \{ y \, | \, d(z_1, y) \le
d(z_2, y) \} \cup \{ y \, | \, d(z_2, y) \le d(z_1, y) \}. \ee For
$y \in D_1$, from (\ref{ddGam}), we have the following bounds. If
$d(z_1, z_2) \le d(z_1, y)$, then
 \be
 \lab{ddGcase1} |\d_{x_j} \d_r
\Gamma(z_1, y)- \d_{x_j} \d_r \Gamma(z_2, y)| \le C [d(z_1,
z_2)/d(z_1, y)]^{-\a_0} d(z_1, y)^{-n}; \ee If $d(z_1, z_2) \ge
d(z_1, y)$, then \be |\d_{x_j} \d_r \Gamma(z_1, y)- \d_{x_j} \d_r
\Gamma(z_2, y)| \le C  d(z_1, y)^{-n}. \ee Similarly

 Therefore, for any positive
number $\a_1<\a_0$, if $d(z_1, z_2) \le d(z_1, y)$, then \be
\lab{ddGcase2}|\d_{x_j} \d_r \Gamma(z_1, y)- \d_{x_j} \d_r
\Gamma(z_2, y)| \le C [d(z_1, z_2)/d(z_1, y)]^{\a_1} d(z_1,
y)^{-n}; \ee If $d(z_1, z_2) \ge d(z_1, y)$, then \be |\d_{x_j}
\d_r \Gamma(z_1, y)- \d_{x_j} \d_r \Gamma(z_2, y)| \le C
\frac{d(z_1, z_2)^{\a_1}}{d(z_1, y)^{n+\a_1}}. \ee Similar, for $y
\in D_1$, inequalities (\ref{ddGcase1}) and (\ref{ddGcase2}) still
hold after switching $z_1$ with $z_2$.

Take $q=z_1$ in (\ref{i4=}). By Lemma \ref{lec1apot} and
(\ref{thetajiehpq}), we also have \be |\p_s u(y) - \p_s u(z_1) |
\le C d^{\a_0}(z_1, y). \ee From (\ref{ddGcase1}) and
(\ref{ddGcase2}), we find that \be \al
&|I_{41}(z_1) - I_{41}(z_2)| d(z_1, z_2)^{-\a_1}\\
&\qquad
 \le C \int_{D_1} \frac{1}{d(z_1, y)^{n-\a_0+\a_1}} dg(y) +
 C \int_{D_2} \frac{1}{d(z_2, y)^{n-\a_0+\a_1}} dg(y)\le C.
 \eal
\ee From this and (\ref{i423ca}), we deduce that $\Vert I_4
\Vert_{C^{\a_1}(B(m, 1/2))} \le C$. Hence we have proven that
(\ref{hic1a0}) is true. This proves the lemma in Case 1.
\medskip

{\it Case 2.} Suppose Condition 2 holds.
\medskip

The proof in this case differs with that of Case 1 only in the paragraph below
(\ref{nolargerrad}).

Assuming the conclusion of the lemma is false, then there exists a sequence of pointed manifolds
$({\M}_i, x_i, g_i)$ which satisfies the basic assumptions, but, for some
$p_i \in B(x_i, r, g_i) \subset {\bf M}_i$, we have
\be
\al
\frac{r^{\theta, \a}_{g_i}(p_i) }{dist( p_i, \, \p B(x_i, r, g_i))}
&=  \inf \{ \frac{r^{\theta, \a}_{g_i}(y) }{dist( y, \, \p B(x_i, r, g_i))}  \, | \,
 y \in B(x_i, r, g_i)
 \} \\
&=  \delta_i \to 0.
\eal
\ee
Again, write $
\lambda_i = dist( p_i, \, \p B(x_i, r, g_i))
$
and consider the scaled metric
$
h_i = (\delta_i \lambda_i)^{-2} g_i.
$
Then , for every point
$y \in B(p_i,  \delta^{-1}_i, h_i)$, the ball $B(y, 1, h_i)$ is contained in a harmonic
coordinate such that  on this same ball and under $h_i$, we have
\be
\lab{hith-alcond2}
  e^{-\theta} I \le ((h_i)_{pq}) \le e^\theta  I  \quad \text{and}
 \quad  \sup_{p, q} ( \Vert (h_i)_{pq} \Vert_{C^{1, \a}}  ) \le e^{\theta}.
\ee However, for any number $\rho>1$,
\be
\lab{nolargerrad2}
\text{there is no harmonic coordinate system
containing} \quad B(p_i, \rho, h_i)
\ee which satisfies (\ref{hith-alcond2}) on $B(p_i, \rho, h_i)$ under the metric $h_i$.

By Lemma 2.1 in [An], there is a subsequence of the triple $\{
B(p_i, 2 \delta^{-1}_i, h_i), p_i, h_i \}$, still denoted by the
same notation, which converges in $C^{1, \a'}_{loc}$ topology to a
complete, pointed manifold $(N, z, h)$. Here $\a'$ is any positive
number strictly less than $\a$. Also $h$ is a $C^{1, \a}$ metric.
Since the original manifold has bounded Ricci potential, we have
shown that the limit manifold is Ricci flat.

In the previous case, we use Condition 1 i.e. the lower bound of
the injectivity radius and Cheeger-Gromoll splitting theorem to
show that the limit manifold $(N, z, h)$ is Euclidean.  Now,
Condition 2 implies that the ratio between the volume of the balls
in the limit manifold and the volume the Euclidean balls with the
same radius is close to $1$. Therefore $(N, z, h)$ must be the
Euclidean space by \cite{An:1}. Afterwards, going exactly as Case
1, we can finish Case 2.
\medskip

{\it Case 3.} Suppose Condition 3 holds.
\medskip

 If the conclusion of the lemma is false,
then again there exists a sequence of pointed manifolds
$({\M}_i, x_i, g_i)$ which satisfies the basic assumptions,  such that
\be
\lab{etaito0}
\sup_{0<\tau \le \tau_0} \mu(g_i, \tau, B(x_i, r, g_i)) = \eta_i \to 0.
\ee But, for some
$p_i \in B(x_i, r, g_i) \subset {\bf M}_i$, we have
\be
\al
\frac{r^{\theta, \a}_{g_i}(p_i) }{dist( p_i, \, \p B(x_i, r, g_i))}
&=  \inf \{ \frac{r^{\theta, \a}_{g_i}(y) }{dist( y, \, \p B(x_i, r, g_i))}  \, | \,
 y \in B(x_i, r, g_i)
 \} \\
&\equiv  \delta_i \to 0.
\eal
\ee

Consider the scaled metric
$
h_i = (\delta_i \lambda_i)^{-2} g_i
$ where $
\lambda_i = dist( p_i, \, \p B(x_i, r, g_i))
$ again.
Then , for every point
$y \in B(p_i,  \delta^{-1}_i, h_i)$, the ball $B(y, 1, h_i)$ is contained in a harmonic
coordinate such that  on this same ball and under $h_i$, we have
\be
\lab{hith-alcond3}
  e^{-\theta} I \le ((h_i)_{pq}) \le e^\theta  I  \quad \text{and}
 \quad  \sup_{p, q} ( \Vert (h_i)_{pq} \Vert_{C^{1, \a}}  ) \le e^{\theta}.
\ee However, for any number $\rho>1$,
\be
\lab{nolargerrad3}
\text{there is no harmonic coordinate system
containing} \quad B(p_i, \rho, h_i)
\ee which satisfies (\ref{hith-alcond3}) on $B(p_i, \rho, h_i)$ under the metric $h_i$.

By Lemma 2.1 in [An], there is a subsequence of the triple
 $\{ B(p_i,  \delta^{-1}_i, h_i), p_i, h_i \}$, still denoted by the same notation,
 which converges in $C^{1, \a'}_{loc}$ topology to a complete,
 pointed manifold $(N, z, h)$. Here $\a'$ is any positive
number strictly less than $\a$. Also $h$ is a $C^{1, \a}$ metric.  Since the original manifold has
 bounded Ricci potential, it is easy to see that the limit manifold is Ricci
flat. It is well known that the following scaling property holds
\be
\al
&\mu(h_i, \frac{1}{2}, B(p_i,  \delta^{-1}_i, h_i) )=
\mu((\delta_i \lambda_i)^2 h_i, \frac{1}{2}(\delta_i \lambda_i)^{2}, B(p_i,
\lambda_i ,g_i))\\
&=\mu(g_i, \frac{1}{2}(\delta_i \lambda_i)^{2}, B(p_i,  \lambda_i ,g_i))
\ge \mu(g_i, \frac{1}{2}(\delta_i \lambda_i)^{2}, B(x_i,  r ,g_i)).
\eal
\ee The last inequality holds since $B(p_i,  \lambda_i ,g_i) \subset B(x_i,  r ,g_i)$
by the definition $\lambda_i=dist(p_i, \p B(x_i, r, g_i))$.

Since $\tau_0$ is a fixed number, if $i$ is sufficiently large,
then $\frac{1}{2}(\delta_i \lambda_i)^{2} \le \tau_0$. By (\ref{etaito0}) we know that
\be
\mu(h_i, \frac{1}{2}, B(p_i,  \delta^{-1}_i, h_i) ) \ge - \eta_i
\ee Since $ \delta^{-1}_i \to \infty$ and $\eta_i \to 0$ when $i \to \infty$, this implies
\be
\mu(h, \frac{1}{2}, N) \ge 0.
\ee
  According to Theorem 16.35 in the book \cite{C++:1}, we know that
 $(N, z, h)$ is in fact
the standard Euclidean space.  We should mention that this result
 was attributed to the paper \cite{BCL:1} Corollary 1.6,
where it was stated slightly differently.

 Now we can finish the proof of the lemma just like
Case 1.
\qed

\section{proof of theorems }

{\it Proof of Theorem \ref{thmconv}}.

With Lemma  \ref{leharmrad} in hand, we know that each manifold in
the class SP has a harmonic coordinate atlas with uniform lower
bound for the radii and uniform $C^{1, \a}$ bound for the metric.
Hence Theorem \ref{thmconv} follows immediately from the general
result Lemma 2.1 in \cite{An:1}. \qed

Next we give

\medskip
{\it Proof of Theorem \ref{thSLocKRF}.}
 Under either one of the
conditions, by Lemma  \ref{leharmrad}, for fixed small $\theta>0$ and each point $p \in \M$,
there is a fixed number $a>0$ such that the ball $B(p, a r, 0)$
 is contained in a harmonic coordinates.
Moreover, the metric $g_{pq}$  satisfies (\ref{gpqC1a}).
  This shows, there exists a smaller positive number $b$ such that the
 isoperimetric constant of the ball $B(p, b r, 0)$ is so close to
 the Euclidean one  that Perelman's pseudolocality theorem \cite{P:1} Section 10
 can be
 applied.  The theorem follows.
\qed

We finish by pointing out a possible smooth convergence result
when the initial metric has entropy close to $0$.

{\it Observation.
 Let $({\M}, g(t))$,  $\partial_t g_{i\g j} = -  R_{i\g j} + g_{i\g{j}}$, be a K\"ahler Ricci flow on a $n$ real
 dimensional compact, K\"ahler manifold with positive first Chern class.
 Suppose also that the initial metric is so scaled that the unnormalized
flow blows up exactly at time $\frac{1}{2}$. There exists a
positive number $\eta_0$ with the following property. If
$\mu(g(0), \frac{1}{2}, {\M}) \ge -\eta_0$, then every sequence
$\{({\M}, g(t_k))\}$, $t_k \to \infty$, sub-converges in
$C^\infty$ topology to a gradient K\"ahler Ricci soliton. Moreover
the curvature tensor is uniformly bounded for all time. }

One may wonder for what initial metrics the unnormalized flow
blows up exactly at time $\frac{1}{2}$?  By Cao's result
\cite{Ca:1} that the normalized flow exists for all time, after a
simple change of time variable, we know that
 this is always true
if $g(0)$ has canonical K\"ahler class, i.e. $2 \pi c_1({\bf M})$,
as its K\"ahler class.

 At this moment, we are not sure if there exist manifolds satisfying the conditions of the
 observation since we do not know the size of $\eta_0$.
Since the argument reveals a connection between the entropies at
initial time and large time, we present it for the interested
reader.

Let $\tilde g$ and $\tilde t$ denote the metric and time of the corresponding un-normalized
Ricci flow, which is made to blow up when $\tilde t = 1/2$. Then we have the relation
\be
t= - \ln (1 - 2 \tilde t), \qquad g(t) = \frac{1}{1- 2 \tilde t} \,  \tilde g ( \tilde t).
\ee From \cite{P:1}, $\forall \e>0$,
\be
\mu(\tilde g ( \tilde t), \e) \ge \mu(g(0), \tilde t + \e),
\ee and therefore
\be
\mu(\frac{1}{1-2 \tilde t} \tilde g ( \tilde t), \frac{1}{1-2 \tilde t} \e) \ge \mu(g(0), \tilde t + \e)
\ee which can be converted to
\be
\mu(g(t), e^t \e) \ge \mu(g(0), \frac{1}{2} - \frac{1}{2} e^{-t} + \e).
\ee Writing $\tau = e^t \e$, then this inequality becomes
\be
\mu(g(t), \tau) \ge \mu(g(0), \frac{1}{2} - \frac{1}{2} e^{-t} + e^{-t} \tau).
\ee  Since $({\M}, g(0))$ is a fixed manifold, given any fixed $\tau_0$, if $\tau \in (0, \tau_0]$,
we can then deduce
\be
\mu(g(t), \tau) \ge \mu(g(0), \frac{1}{2}) + o(1), \qquad t \to \infty.
\ee  By our assumption that $\mu(g(0), \frac{1}{2}) \ge -\eta$, with $\eta$ small, we know that
the family $\{ {\M}, g(t) \}$, $t$ large, satisfies the conditions of Theorem \ref{thmconv} (b). This
proves  sub-convergence in $C^{1, \a}$ topology.

Now we prove the convergence in $C^\infty$ topology. Pick a sequence
of times $t_k$ which go to $0$ when $k \to \infty$. Let $\delta>0$ be a small
number to be chosen later. Then a subsequence of $\{(M, g(t_k-\delta))\}$,
still identified by the same notation, converges in $C^{1, \a}$ topology.
By Lemma \ref{leharmrad}, for fixed $\theta>0$ and each point $p \in \M$,
there is a fixed number $r_0$ such that the ball $B(p, r_0, t_k-\delta)$
 is contained in a harmonic coordinates.
Moreover, the metric $g_{pq}$ in each of the ball satisfies, \be
\lab{gpqC1a} a). \quad  e^{-\theta} I \le (g_{pq}) \le e^\theta  I
; \quad b). \quad  \sup_{p, q} ( \Vert g_{pq} \Vert_{C^{1}} \, r_0
) \le e^{\theta}. \qquad c). \sup_{p, q} ( \Vert g_{pq}
\Vert_{C^{1, \a}} r^{1 + \a}_0 ) \le e^{\theta}. \ee This shows,
by choosing $\theta$ small, there exists a small positive number
$\e$ such that the
 isoperimetric constant of the ball $B(p, \e r_0, g(t_k-\delta))$ is so close to
 the Euclidean one  that Perelman's pseudolocality theorem \cite{P:1} can be
 applied. This implies, if one chooses $\delta=c_0 (\e r_0)^2$ with $c_0>0$
sufficiently small, then the curvature tensor satisfies
\be
\lab{Rmjie}
|Rm(p, t)| \le C (\e r_0)^{-2}, \qquad t \in [t_k-\frac{\delta}{2}, t_k+
\frac{\delta}{2}].
\ee Since $p$ is arbitrary, we find that $|Rm(\cdot, t)|$,
$t \in [t_k-\frac{\delta}{2}, t_k+\frac{\delta}{2}]$,  is uniformly bounded.
Hence the convergence is actually in $C^\infty$ topology. Using
the result in \cite{Se:1}, we know that the limit manifold is a gradient
K\"ahler Ricci soliton.  See also Lemma 3.5 in \cite{TZhu2:1}  for a detailed proof.

Finally, since $\{ t_k \}$ is an arbitrary time sequence, from
(\ref{Rmjie}), we see that the curvature tensor $Rm$ is uniformly
bounded for all time. This confirms the observation.

\medskip

{\bf Acknowledgment.} Q. S. Z. would like to thank Professors X. X. Chen,
G. F. Wei and  Zhenlei Zhang for
helpful conversations.

\bigskip

\noindent e-mail:  tian@math.princeton.edu and qizhang@math.ucr.edu

\enddocument